\documentclass[journal]{IEEEtran}
\usepackage{amsfonts}
\usepackage{bm}
\usepackage{amsmath,amsfonts,amsthm,amssymb}
\usepackage{setspace}
\usepackage{Tabbing}
\usepackage{fancyhdr}
\usepackage{lastpage}
\usepackage{extramarks}
\usepackage{chngpage}
\usepackage{algorithmic}
\usepackage{soul,color}
\usepackage{epsfig}
\usepackage{graphicx,float,wrapfig,subfigure}
\usepackage{dsfont}
\usepackage{longtable}
\usepackage{wrapfig}
\usepackage{stfloats}
\usepackage{cite}
\usepackage{graphicx}
\usepackage{array}
\usepackage{multirow}{\tiny }
\usepackage{multicol}
\usepackage{colortbl}
\usepackage{tabularx}
\usepackage{mdwmath}
\usepackage{mdwtab}
\usepackage{color}
\usepackage{verbatim}
\usepackage{amsmath}
\usepackage{tikz}
\usetikzlibrary{calc}
\usepackage{flowchart}
\usetikzlibrary{shapes.geometric, arrows}
\usepackage{overpic}
\usepackage{epstopdf}
\usepackage{url}
\usepackage{bbm}
\usepackage[colorlinks,linkcolor=black,anchorcolor=black,citecolor=black,urlcolor=black]{hyperref} 

\usepackage{makecell}
\usepackage{textcomp}
\usepackage{algorithm}
\usepackage{algorithmic}

\usepackage{threeparttable}

\newtheorem{proposition}{Proposition}

\newtheorem{remark}{Remark}
\DeclareMathOperator*{\argmin}{argmin}

\usepackage{color}
\makeatletter %
\let\myorg@bibitem\bibitem
\def\bibitem#1#2\par{%
	\@ifundefined{bibitem@#1}{%
		\myorg@bibitem{#1}#2\par
	}{%
		\begingroup
		\color{\csname bibitem@#1\endcsname}%
		\myorg@bibitem{#1}#2\par
		\endgroup
	}%
}

\makeatother %

\allowdisplaybreaks

\usepackage{footnote}
\makesavenoteenv{tabular}
\makesavenoteenv{table}

\begin{document}

\title{
Deep-quantile-regression-based surrogate model for joint chance-constrained optimal power flow with renewable generation
	}

\author{
Ge~Chen,~\IEEEmembership{Graduate Student Member,~IEEE,}
Hongcai~Zhang,~\IEEEmembership{Member,~IEEE,}
Hongxun~Hui,~\IEEEmembership{Member,~IEEE,}
and~Yonghua~Song,~\IEEEmembership{Fellow,~IEEE}
\vspace{-8mm}

}

\maketitle

\begin{abstract}
Joint chance-constrained optimal power flow (JCC-OPF) is a promising tool to manage uncertainties from distributed renewable generation. However, most existing works are based on power flow equations, which require accurate network parameters that may be unobservable in many distribution systems. To address this issue, this paper proposes a learning-based surrogate model for JCC-OPF with renewable generation. This model equivalently converts joint chance constraints in quantile-based forms and introduces deep quantile regression to replicate them, in which a multi-layer perceptron (MLP) is trained with a special loss function to predict the quantile of constraint violations. Another MLP is trained to predict the expected power loss. Then, the JCC-OPF can be formulated without network parameters by reformulating these two MLPs into mixed-integer linear constraints. To further improve its performance, two pre-processing steps, i.e., data augmentation and calibration, are developed. The former trains a simulator to generate more training samples for enhancing the prediction accuracy of MLPs. The latter designs a positive parameter to calibrate the predictions of MLPs so that the feasibility of solutions can be guaranteed. Numerical experiments based on the IEEE 33- and 123-bus systems validate that the proposed model can achieve desirable feasibility and optimality simultaneously with no need for network parameters.
\end{abstract}
\begin{IEEEkeywords}
Optimal power flow, joint chance constraints, deep quantile regression, distribution network, distributed renewable generation.
\end{IEEEkeywords}

\section{Introduction} \label{sec_intro}
\IEEEPARstart{O}{ptimal} power flow (OPF) plays a critical role in the operation of distribution networks \cite{abdi2017review}. By solving OPF, network operators can find the most economical dispatch strategy while ensuring operational security. However, since the distributed generators (DGs), such as wind turbines and PV plants, have been increasingly integrated into distribution networks \cite{9709098}, considerable uncertainties are introduced in OPF, which dramatically increases the difficulty for solving OPF \cite{zakaria2020uncertainty}.

Chance constrained programming (CCP) is a promising method to account for the uncertainties from DGs in OPF \cite{geng2019data}. It allows constraint violation with a small probability so that operators can effectively balance robustness and optimality based on their preferences. Many recent efforts have also been made to use CCP to describe the impacts of uncertainties in OPF. References \cite{8515118,8688432} combined CCP to the linearized DistFlow model to coordinate uncertain DGs with flexible resources in distribution networks. Reference \cite{8600344} applied CCP to a linearized AC OPF model to schedule wind generation. Nevertheless, the above papers used individual chance constraints to control the violation probability of every critical constraint. This individual manner may not guarantee the joint satisfaction probability of all critical constraints, while this joint probability is more concerned for operators because they need to ensure the security of entire systems \cite{7973099}. 

Conversely, the joint chance-constrained optimal power flow (JCC-OPF) directly restricts the joint satisfaction probability of all critical constraints \cite{9122389}. Hence, it is preferable to ensure the system-level security and attract increasing attention in recent years. References \cite{8355588,8662704} employed JCC-OPF to restrict the joint probability of critical constraint violations, where Bonferroni approximation was used to convert the intractable joint chance constraints into solvable individual ones. Reference \cite{8528856} proposed a joint chance-constrained linearized DistFlow model to schedule the reactive power compensation for distribution networks. References \cite{8060613, 8626040} combined JCC-OPF with a scenario approach to approximately reformulate the probability constraints into tractable deterministic forms.  Nevertheless, most existing works still face two challenges:
\begin{enumerate}
\item Most published works, including \cite{9122389,8355588,8662704,8528856,8060613, 8626040}, are based on power flow equations, which require accurate power network parameters including topology of network and impedance of each line etc. However, 
these parameters are often unavailable in many distribution networks due to the unaware topology changes or inaccurate data maintenance \cite{7875102}.
\item The impacts of uncertainties on bus voltages and branch currents are hard to quantify because the OPF model is non-convex. Existing papers usually introduce approximations (e.g. the linearized DistFlow used in \cite{8528856}) or relaxations (e.g. the semi-definite relaxation used in \cite{8060613, 8626040}) so that the impacts of uncertainties are convenient to describe. However, these approximation or relaxation models may lead to overly conservative or even infeasible solutions. For instance, in a radial distribution network with high DG penetration, reverse power flows may occur. In that case, the semi-definite relaxation, which is equivalent to a second-order cone (SOCP) relaxation in radial networks, is not exact and may not ensure feasibility \cite{6815671}.
\end{enumerate}

Since collecting historical data (e.g. power injections, bus voltages, and branch currents) becomes easier and cheaper nowadays, learning-based methods may be a potential choice to bypass the above two challenges because they can train tractable surrogate models without the network parameters to replace the non-convex power flow model \cite{9265482}. Generally speaking, the existing learning-based methods can be divided into the following three categories.

\subsubsection{Learn optimal solutions}
Methods in this category usually train neural networks to directly learn the optimal solution of OPF. 
For example, reference \cite{9205647} trained a neural network to build a mapping from power demands to the optimal solution of DC OPF. In references \cite{chatzos2020high,9335481}, this learning-based method was combined with the Lagrangian dual approach to improve the feasibility of solutions. References \cite{8810819,9599403} trained neural networks to learn active constraints of OPF. Then, optimal solutions were obtained by solving the equations formed by these active constraints. Generally speaking, these methods could dramatically reduce the solving time of OPF. However, they need optimal solutions of OPF as the training labels, so the network parameters are still indispensable.

\subsubsection{Learn feasibility conditions}
Methods in this category usually learn to replicate OPF constraints by training neural networks. For instance, in references \cite{9302963,venzke2020neural}, binary classifiers were trained to judge whether a given strategy can satisfy all constraints or not. Then, the trained classifiers were equivalently reformulated as mixed-integer linear constraints so that the OPF problem can be replicated with no need for building any power flow model. Reference \cite{9502573} replaced the binary classifiers with a regression neural network that can predict the maximum constraint violation to improve the feasibility of solutions. The above methods only require solutions of power flow equations instead of optimal solutions of OPF as training labels, so the requirement of network parameters can be bypassed. Moreover, desirable optimality can be also achieved since the mixed-integer linear replication of OPF can be efficiently solved by the Branch-and-Bound algorithm.
However, it is difficult for these methods to quantitatively evaluate the impacts of uncertainties from DGs (inputs of neural networks) on the constraint violations (outputs of neural networks) because of the nonlinear activation functions in neural networks. Thus, they are not applicable to JCC-OPF.

\subsubsection{Reinforcement learning}
Reinforcement learning (RL) trains agents how to act in a specific environment to maximize the cumulative reward (e.g. the opposite of energy purchasing from upper-level grids). In reference \cite{9069289}, RL was applied to solve OPF, while the Lagrangian dual approach was combined to improve feasibility. In reference \cite{9275611}, behavior cloning was combined with RL to generate a desirable initial start so that the training process of agents can be accelerated. However, RL needs ``trial-and-errors" to train agents, which may be unacceptable in practical operations of distribution networks. 

Moreover, it is also challenging for all the learning-based methods above to handle joint chance constraints. If they want to learn the characteristics of joint chance constraints and train surrogate models to replace them, then their training sets must contain enough samples of statistical results (e.g., the quantile of constraint violations). However, these samples are difficult to collect because only realizations of uncertainties instead of the quantile can be observed in practice. In fact, to the best of our knowledge, none of these learning-based methods have been successfully extended to the JCC-OPF problem. 

To overcome the aforementioned challenges, this paper proposes a novel learning-based surrogate model for JCC-OPF with renewable generation. Two pre-processing steps, i.e., data augmentation and calibration, are further designed to improve the performance of the proposed model. The specific contributions are threefold:
\begin{enumerate}
    \item We propose a learning-based surrogate model for JCC-OPF. The proposed model first re-expresses joint chance constraints in quantile-based forms and introduces deep quantile regression to replicate them, where a multi-layer perceptron (MLP) is trained based on a special loss function to predict the quantile of the maximum constraint violation. Then, another MLP is trained based on mean squared errors to predict the expected power loss. By reformulating these two MLPs into mixed-integer linear constraints, the proposed surrogate model can be established. Since the surrogate model only requires historical data to train MLPs but does not need to build exact power flow models, the requirement of network parameters can be bypassed. Moreover, this model can be efficiently solved by the Branch-and-Bound algorithm with guaranteed optimality.
    \item Considering that the historical dataset may not contain enough training samples to reflect the true distribution of constraint violations, the prediction accuracy of the quantile regression may be undesirable. To address this issue, a data augmentation step is designed. This step uses the historical dataset to train a regressor as a simulator based on XGBoost. With this simulator, more training samples can be generated for the previous MLPs to improve the accuracy of the quantile regression.
    \item Since the deep quantile regression may have prediction errors and harm the feasibility of solutions, a calibration step is further designed.
    In this step, we first demonstrate that the underestimation for the quantile of constraint violations may lead to infeasible solutions. Then, a positive constant, i.e., calibration parameter, is designed to calibrate the outputs of the deep quantile regression to avoid the harmful underestimation so that the feasibility of solutions can be improved.
\end{enumerate}

The remaining parts are organized as follows. Section \ref{sec_formulation} describes the formulation of the JCC-OPF problem. Section \ref{sec_solution} introduces the proposed learning-based surrogate model in detail. Section \ref{sec_case} demonstrates simulation results and Section \ref{sec_conclusion} concludes this paper.

\section{Formulation of JCC-OPF} \label{sec_formulation}
This paper develops a learning-based surrogate model of the JCC-OPF problem for a distribution network without network parameters. In this section, we first present the detailed formulation of the JCC-OPF problem.

\subsubsection{Power injections}
By using $i \in \mathcal{V}$ to index buses, the active and reactive power injections on each bus, i.e., $\bm p \in \mathbb{R}^{|\mathcal{V}|}$ and $\bm q \in \mathbb{R}^{|\mathcal{V}|}$, can be expressed as:
\begin{align}
\bm p = - \bm p^\text{d} + \bm p^\text{DG}, \quad \bm q =  - \bm q^\text{d} + \bm \phi * \bm p^\text{DG}, \label{eqn_injection}
\end{align}
where $\bm p^\text{d}$ and $\bm q^\text{d}$ represent the active and reactive power demands on each bus. Variable $\bm p^\text{DG}$ is the actual used active power from DGs. Parameter $\bm \phi$ is a ratio of the actual active power of DG to its reactive power. Operator $*$ denotes the element-wise multiplication. The actual used active power $\bm p^\text{DG}$ can be expressed by:
\begin{align} 
	\bm p^\text{DG} = \bm \lambda * \bm G^\text{DG},
\end{align}
where $\bm \lambda$ and $G^\text{DG}$ are the actual utilization rate and maximum available value of DG. In practice, the value of $\bm G^\text{DG}$ is uncertain, which can be expressed as follows:
\begin{align} 
	\bm G^\text{DG} = \overline{\bm G}^\text{DG} * (\bm 1 + \bm \omega),
\end{align}
where $\overline{\bm G}^\text{DG}$ represents the nominal available DG obtained by predictions and $\bm \omega$ is the corresponding uncertain level. 


\subsubsection{Power flow model}
The power flow model of a radial network can be expressed by DistFlow \cite{19266}, as follows:
\begin{align}
\begin{cases}
\sum_{k \in \mathcal{C}_{j}} P_{jk} = p_{j} + P_{ij} - r_{ij}I_{ij}^2, \\
\sum_{k \in \mathcal{C}_{j}} Q_{jk} = q_{j} + Q_{ij} - x_{ij}I_{ij}^2, \\
V_{j}^2=V_{i}^2-2(r_{ij}P_{ij}+x_{ij}Q_{ij})\\
\quad \quad \quad \quad \quad+ (r_{ij}^2 + x_{ij}^2)I_{ij}^2, \\
I_{ij}^2 = \frac{P_{ij}^2 + Q_{ij}^2}{V_{i}^2},\\
\end{cases} \forall (i,j) \in \mathcal{B}, \label{eqn_distflow}
\end{align}
where $P_{ij}$ and $Q_{ij}$ are the active and reactive power flows on branch $(i, j)$, respectively; $V_{i}$ and $I_{ij}$ are the magnitudes of the voltage at bus $i$ and current on branch $(i, j)$, respectively; $r_{ij}$ and $x_{ij}$ denotes the resistance and reactance of branch $(i,j)$, respectively. Set $\mathcal{C}_{j}$ contains the child bus indexes of bus $j$. Set $\mathcal{B}$ represents the index set of branches in this network.

\subsubsection{Security constraints}
To ensure operation security, the magnitudes of all bus voltages and branch currents shall maintain in the corresponding allowable ranges.
According to (\ref{eqn_injection})-(\ref{eqn_distflow}), the uncertainties from DGs also affect the bus voltages and branch currents. To better balance optimality and feasibility of OPF solutions, a joint chance constraint is employed to describe the voltage and current limitations:
\begin{align}
\mathbb{P}_{\bm \omega}\left( \bm V_\text{min} \leq \bm V \leq \bm V_\text{max}, \ \bm I\leq \bm I_\text{max} \right) \geq 1 - \epsilon, \label{eqn_JCC}
\end{align}
where $\bm V$ and $\bm I$ are the vector forms of $V_{i}$ and $I_{ij}$; $\epsilon$ is the risk parameter. Note here we use a joint chance constraint instead of individual ones because this joint manner can better guarantee the security of the entire system \cite{9122389}.

\subsubsection{Energy purchasing}
The energy purchasing from the upper-level grid $G$ is equal to the net power at the substation and can be calculated based on the power balance of a distribution network:
\begin{align}
G  = \bm 1^\intercal \bm p^\text{d} + p^\text{loss} - \bm 1^\intercal (\bm \lambda * \bm G^\text{DG}), \label{eqn_balance}
\end{align}
where $p^\text{loss}$ is the total power loss and can be calculated by:
\begin{align}
p^\text{loss} = \sum_{(i,j) \in \mathcal{B}} r_{ij} I_{ij}^2. \label{eqn_loss}
\end{align}
Finally, the JCC-OPF is formulated as:
\begin{align} 
	&\min_{\bm \lambda, G} \quad \mathbb{E}_{\bm \omega}(G),\quad \text{s.t.:} \text{ Eqs. (\ref{eqn_injection})-(\ref{eqn_loss})}. \tag{$\textbf{P1}$}
\end{align}

\section{Solution Methodology} \label{sec_solution}
As mentioned in Section \ref{sec_intro}, formulating \textbf{P1} can be challenging because the network parameters are often unavailable. Fortunately, collecting historical operation data becomes easier and cheaper due to the widespread use of smart meters. Therefore, we propose a learning-based surrogate model to address the aforementioned challenge. This model first introduces deep quantile regression to replicate the intractable joint chance constraints. Meanwhile, another neural network is trained to predict the expected power loss. Then, by reformulating the trained neural networks into mixed-integer linear constraints, \textbf{P1} can be replicated. Since the proposed model only requires historical data to train neural networks but does not need to build power flow models, the requirement of the network parameters can be bypassed. 

\subsection{Deep quantile regression to replicate chance constraints}
\subsubsection{Motivation}
The joint chance constraint (\ref{eqn_JCC}) can be equivalently reformulated into quantile-based deterministic forms to eliminate the intractable probability operator. 
Specifically, we define $\bm x$ as the nominal active and reactive power injections on each bus (except the slack bus): 
\begin{align} 
    \bm x = [\bm p, \  \bm q]. \label{eqn_x_define}
\end{align}
We also define a new variable $h$ to denote the maximum violation of the OPF constraints:
\begin{align} 
	h(\bm x, \bm \omega) = \max \{\bm V_\text{min} - \bm V, \bm V - \bm V_\text{max}, \bm I - \bm I_\text{max}\}. \label{eqn_h_define}
\end{align}
Note the impacts of the uncertainty from DGs, i.e., $\bm \omega$, has been implied in the samples of $h$ because both the voltage $\bm V$ and current $\bm I$ are affected by $\bm \omega$. Based on (\ref{eqn_h_define}), the joint chance constraint  (\ref{eqn_JCC}) can be expressed as:
\begin{align}
\mathbb{P}_{\bm \omega}\left( h(\bm x, \bm \omega) \leq 0 \right) \geq 1 - \epsilon, \label{eqn_JCC2}
\end{align}
which can be further equivalently reformulated into the following quantile-based form:
\begin{align}
\mathcal{Q}_{\bm \omega}^{1-\epsilon}(h(\bm x, \bm \omega)) \leq 0, \label{eqn_quantile_JCC}
\end{align}
where $\mathcal{Q}_{\bm \omega}^{1-\epsilon}(h(\bm x, \bm \omega))$ is the $1-\epsilon$ quantile of $h$ at a given $\bm x$:
\begin{align}
\mathcal{Q}_{\bm \omega}^{1-\epsilon}(h(\bm x, \bm \omega)) = \inf\{y: \mathbb{P}_{\bm \omega}\left(h(\bm x, \bm \omega)\leq y \right)\geq 1-\epsilon\}. \label{eqn_quantile_definition}
\end{align}
According to (\ref{eqn_quantile_JCC}), if the mapping from $\bm x$ to $\mathcal{Q}_{\bm \omega}^{1-\epsilon}(h(\bm x, \bm \omega))$ can be accurately described by simple relations (e.g. linear functions), then the intractability of the joint chance constraint can be overcome. This motivates us to introduce a powerful deep learning technique, deep quantile regression, to predict the quantile of constraint violations. 

\subsubsection{Introduction of deep quantile regression}
Traditional regression is a process to model the relationship between dependent output and independent variables. For example, based on the dataset $\{(\bm x_n, \bm \omega_n, h_n)\}_{n \in \mathcal{N}}$ ($\mathcal{N}$ is the index set of samples), we can train a regression model $\hat h(\bm x, \bm \omega)$ to predict $h$ with a given $\bm x$ and $\bm \omega$.
The mean squared error is usually used as the loss function in traditional regression models:
\begin{align}
\text{Loss}^\text{R} = (h - \hat{h}(\bm x, \bm \omega))^2.
\end{align}
However, it is hard for traditional regression models to accurately predict the quantile $\mathcal{Q}_{\bm \omega}^{1-\epsilon}(h(\bm x, \bm \omega))$. This is because they need sufficient samples of the quantile as the training labels, which are usually difficult to collect. In practice, for a specific $\bm x=\bm x_n$, we may only observe one realization of $h$, i.e., $h_n=h(\bm x_n, \bm \omega_n)$, instead of the quantile at $\bm x_n$ (other samples usually have different $\bm x$). Without enough training labels, traditional regression models can not work well. 

Conversely, deep quantile regression can directly predict the quantile based on realizations of uncertainties, i.e., $h_n$, instead of quantile samples \cite{hao2007quantile}. This advantage results from a specially designed loss function, as follows:
\begin{align}
\text{Loss}^\text{QR} = \mu \cdot(1-\epsilon-\mathbb{I}(\mu\leq 0)). \label{eqn_loss_quantile}
\end{align}
Here $\mu=h-\mathcal{\hat Q}^{1-\epsilon} (\bm x)$, where $\mathcal{\hat Q}^{1-\epsilon} (\bm x)$ is the prediction of $\mathcal{Q}_{\bm \omega}^{1-\epsilon}(h(\bm x, \bm \omega))$ given by the quantile regression. Symbol $\mathbb{I}(\cdot)$ denotes the indicator function. The following \textbf{Proposition} proves that we can predict the quantile without the samples of $\mathcal{Q}_{\bm \omega}^{1-\epsilon}(h(\bm x, \bm \omega))$ based on the loss function (\ref{eqn_loss_quantile}).

\begin{proposition} \label{proposition_1}
The quantile $\mathcal{Q}_{\bm \omega}^{1-\epsilon}(h(\bm x, \bm \omega))$ can be obtained by minimizing the expectation of (\ref{eqn_loss_quantile}), as follows \cite{hao2007quantile}:
\begin{align}
\mathcal{Q}_{\bm \omega}^{1-\epsilon}(h(\bm x, \bm \omega)) = \argmin_{\mathcal{\hat Q}^{1-\epsilon} (\bm x)} \mathbb E (\rm{Loss}^{\rm{QR}}). \label{eqn_proposition_1}
\end{align}
\end{proposition}
\emph{Proof}: See Appendix \ref{app_1}. 

Based on \textbf{Proposition} \ref{proposition_1}, we can train a quantile regression neural network based on the historical dataset $\{(\bm x_n, h_n)\}_{n \in \mathcal{N}}$ to represent the mapping from $\bm x$ to the target quantile $\mathcal{Q}_{\bm \omega}^{1-\epsilon}(h(\bm x, \bm \omega))$. Note the label $h_n$ in the historical dataset is noisy due to the impacts of uncertainty $\bm \omega$.

\subsubsection{Replication of joint chance constraints}
\begin{figure}
	\centering
	 			\vspace{-4mm}
	\includegraphics[width=0.9\columnwidth]{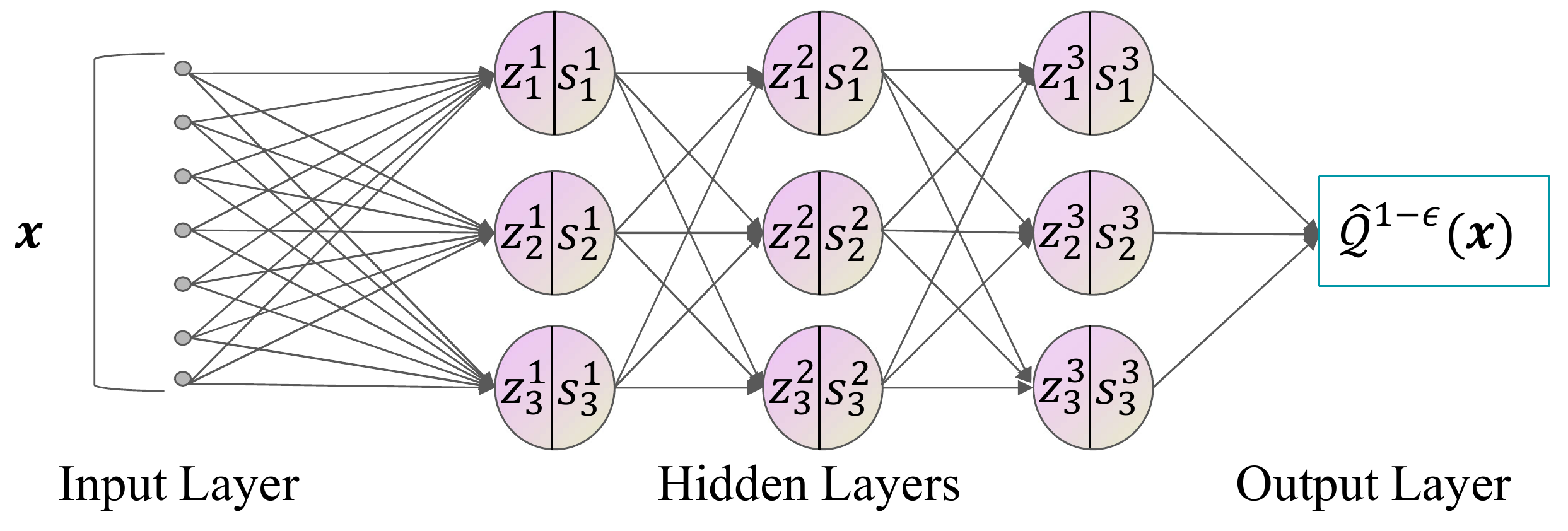}\vspace{-2mm}
	\caption{Structure of an example MLP with 3 hidden layers, where $z_n^l$ and $h_n^l$ denote the outputs of the linear mapping and ReLU of the $n$-th neroun in layer $l$, respectively.
	}
	\label{fig_MLP}
	 		\vspace{-4mm}
\end{figure}

With the noisy dataset $\{(\bm x_n, h_n)\}_{n \in \mathcal{N}}$, we can train a quantile regression network $\mathcal{\hat Q}^{1-\epsilon} (\bm x)$ to predict the target quantile $\mathcal{Q}_{\bm \omega}^{1-\epsilon}(h(\bm x, \bm \omega))$. Then, Eq. (\ref{eqn_quantile_JCC}) is replaced by:
\begin{align}
\mathcal{\hat Q}^{1-\epsilon} (\bm x) \leq 0. \label{eqn_quantile_JCC2}
\end{align}
In this paper, a MLP with ReLU activation functions is chosen as the quantile regression model. For convenience, this MLP is called ``quantile-MLP". A typical MLP is composed of one input layer, $|\mathcal{L}|$ hidden layers, and one output layer, as shown in Fig. \ref{fig_MLP}. Each neuron is made up of a linear mapping and a nonlinear ReLU function. By using $l$ index the hidden layers ($l \in \mathcal{L}$), the target quantile can be estimated by the forward propagation of the trained quantile-MLP:
\begin{align}
&\bm s^{0} = \bm x, \label{eqn_input}\\
\begin{split}
&\bm z^{l} = \bm W^{l} \bm s^{l-1} + \bm b^{l}, 
\forall l \in \mathcal{L}, \label{eqn_z}
\end{split} \\
&\bm s^{l} = \max(\bm z^{l}, 0), \quad \forall l \in \mathcal{L}, \label{eqn_h} \\
&\hat{\mathcal{Q}}^{1-\epsilon}(\bm x) = \bm W^{|\mathcal{L}|+1} \bm s^{|\mathcal{L}|} + \bm b^{|\mathcal{L}|+1}. \label{eqn_output} 
\end{align}
Eq. (\ref{eqn_input}) defines the input layer; Eqs. (\ref{eqn_z}) and (\ref{eqn_h}) represent the linear mapping and ReLU in hidden layers, respectively; Eq. (\ref{eqn_output}) defines the output layer. Vector $\bm z^{l}$ and $\bm s^{l}$ are the outputs of the linear mapping and activation function in hidden layer $l$; matrix $\bm W^l$ and vector $\bm b^{l}$ are the weights and bias of layer $l$, which are parameters to be learned; $\hat{\mathcal{Q}}^{1-\epsilon}(\bm x)$ is the estimation of $\mathcal{Q}_{\bm \omega}^{1-\epsilon}(h(\bm x, \bm \omega))$ given by the quantile-MLP.

\begin{remark}
We can also let the quantile-MLP output multiple quantile values at once to improve its practicability. To realize this, we only need to change its forward propagation into:
\begin{align}
\begin{cases}
\text{Eqs. (\ref{eqn_input})-(\ref{eqn_h})},\\
\left[\hat{\mathcal{Q}}^{1-\epsilon_i}(\bm x), \forall i \in \mathcal{I}\right]^{\intercal} = \bm W^{|\mathcal{L}|+1} \bm s^{|\mathcal{L}|} + \bm b^{|\mathcal{L}|+1}, \label{eqn_output2}
\end{cases}
\end{align}
where $\mathcal{I}$ is the index set of risk parameters. Then, even if multiple quantile values are required, only one single MLP needs to be trained.
\end{remark}

\subsection{Power loss calculation}
Based on (\ref{eqn_balance}), the objective of \textbf{P1} can be calculated by:
\begin{align}
\mathbb E_{\bm \omega} (G) &= \mathbb E_{\bm \omega} \left(\bm 1^\intercal \bm p^\text{d} + p^\text{loss} -  \bm 1^\intercal (\bm \lambda * \bm G^\text{DG})\right), \notag \\ 
& = \bm 1^\intercal \bm p^\text{d} + \mathbb E_{\bm \omega} (p^\text{loss}) - \bm 1^\intercal (\bm \lambda * \overline{\bm G}^\text{DG}). \label{eqn_balance2}
\end{align}
Eq. (\ref{eqn_balance2}) indicates that the expected power loss is required for the calculation of \textbf{P1}'s objective. According to (\ref{eqn_loss}), the power loss calculation requires the magnitudes of branch currents. Thus, the power flow model is still required. To bypass this requirement, another MLP (we call it ``loss-MLP"), is trained to predict the expected power loss.
The samples of nominal power injection $\bm x$ and power loss $p^\text{loss}$ are treated as the features and noisy labels ($p^\text{loss}$ is affected by the uncertainty $\bm \omega$), respectively. The mean squared error is employed as the loss function, as follows:
\begin{align}
\text{Loss}^\text{pl} = (p^\text{loss} - \hat{p}^\text{loss}(\bm x))^2, \label{eqn_loss_powerLoss}
\end{align}
where $\hat{p}^\text{loss}(\bm x)$ is the prediction given by the loss-MLP.
\begin{proposition} \label{proposition_2}
The expected power loss can be obtained by minimizing the expectation of (\ref{eqn_loss_powerLoss}), as follows:
\begin{align}
\mathbb E_{\bm \omega} (p^\text{loss}) = \argmin_{\hat{p}^\text{loss}(\bm x)} \mathbb E (\rm{Loss}^{\rm{pl}}). \label{eqn_proposition_22}
\end{align}
\end{proposition}
\emph{Proof}: See Appendix \ref{app_2}.

After training, the expected power loss can be also predicted based on the forward propagation of the loss-MLP:
\begin{align}
&\bm s_{pl}^{0} = \bm x, \label{eqn_input_loss}\\
\begin{split}
&\bm z_{pl}^{l} = \bm W_{pl}^{l} \bm s_{pl}^{l-1} + \bm b_{pl}^{l}, 
\forall l \in \mathcal{L}_{pl}, \label{eqn_z_loss}
\end{split} \\
&\bm s_{pl}^{l} = \max(\bm z_{pl}^{l}, 0), \quad \forall l \in \mathcal{L}_{pl}, \label{eqn_h_loss} \\
&\hat p^\text{loss} (\bm x) = \bm W_{pl}^{|\mathcal{L}_{pl}|+1} \bm s^{|\mathcal{L}_{pl}|} + \bm b_{pl}^{|\mathcal{L}_{pl}|+1}, \label{eqn_output_loss} 
\end{align}
where the subscript $pl$ is used to mark those variables belonging to the loss-MLP.

\subsection{Tractable reformulation of MLPs}
Once the two MLPs are trained, the quantile of the maximum constraint violation and expected power loss can be predicted by (\ref{eqn_input})-(\ref{eqn_output}) and (\ref{eqn_input_loss})-(\ref{eqn_output_loss}) with no need for building any power flow model. However, Eqs. (\ref{eqn_h}) and (\ref{eqn_h_loss}) are intractable for off-the-shelf solvers due to the maximum operator. To address this, the Big-M reformulation used in \cite{9302963,venzke2020neural,9502573} is employed to convert these intractable constraints into mixed-integer linear forms. Specifically, by introducing auxiliary variables $\bm r^{l}$ and $\bm \mu^{l}$ for each hidden layer, Eqs. (\ref{eqn_z})-(\ref{eqn_h}) can be reformulated as:
\begin{align}
	&\left\{
	\begin{aligned}
	&\bm s^{l} - \bm r^{l}=\bm W^{l} \bm s^{l-1} + \bm b^{l}, \\
	&0 \leq \bm s^{l} \leq M \cdot \bm \mu^{l},\\
	&0 \leq \bm r^{l} \leq M\cdot(1-\bm \mu^{l}),\\
	&\bm \mu^{l} \in \{0,1\}^{N_l},
	 \end{aligned}\right.
	 \label{eqn_reformulation}
\end{align}
where $N_l$ denotes the neuron number in the $l$-th hidden layer of the quantile-MLP. Similarly, Eqs. (\ref{eqn_z_loss})-(\ref{eqn_h_loss}) can be equivalently converted into the same form of (\ref{eqn_reformulation}), which is recorded as follows for convenience:
\begin{align}
	\text{\{Eq. (\ref{eqn_reformulation})\}}_{pl}. \label{eqn_reformulation_loss}
\end{align}
Then, the JCC-OPF problem, i.e., \textbf{P1}, can be replicated by the following learning-based surrogate model:
\begin{align} 
	&\min_{\bm \lambda, G} \quad \mathbb E_{\bm \omega} (G) \tag{$\textbf{P2}$},\\
	&\begin{array}{r@{\quad}r@{}l@{\quad}l}
		\text{s.t.} &&\text{Eqs. (\ref{eqn_x_define}), (\ref{eqn_quantile_JCC2}), (\ref{eqn_input}), (\ref{eqn_output}), (\ref{eqn_balance2}), (\ref{eqn_input_loss}), and (\ref{eqn_output_loss})-(\ref{eqn_reformulation_loss}).} 
	\end{array} \notag
\end{align}
The number of the auxiliary binary variables in \textbf{P2} are the same as the total neuron numbers of the two MLPs. According to our test, with only a few neurons, the two MLPs can already achieve desirable prediction accuracy. Thus, the computational performance of \textbf{P2} is acceptable. This will be also verified by the simulations in Section \ref{sec_sensitive}. 

\begin{remark}
The proposed surrogate model only requires historical samples to train the quantile-MLP and loss-MLP but does not need to build exact power flow model. Thus, the requirement of the network parameters can be bypassed. 
\end{remark}

\subsection{Pre-processing to improve performance}
\subsubsection{Motivation}
If the quantile-MLP is directly trained based on the historical dataset $\{\bm x_n, h_n\}_{n \in \mathcal{N}}$, its prediction accuracy may be poor because the historical dataset may not have enough samples to precisely reflect the true distribution of $h$ at a given $\bm x$. For example, we may only find one sample $\{\bm x_n, h_n\}$ at $\bm x = \bm x_n$ (other samples usually have different $\bm x$). With only one sample, it is hard to learn the true distribution of $h$ at $\bm x = \bm x_n$. Moreover, even if we have an ideal training set to train the quantile-MLP, prediction errors are still inevitable, which may harm the feasibility of the proposed model. To overcome the above challenges, two pre-processing steps, i.e., data augmentation and calibration, are designed.

\subsubsection{Data augmentation}
We design a data augmentation step to construct an ideal training set for the quantile-MLP. Its key idea is very simple: train a simulator based on the historical data and use this simulator to generate more samples as the training set. The detailed procedure is summarized in Table \ref{tab_algorithm1}. Then, at a given $\bm x= \bm x^{(k)}$, multiple labels, i.e., $\{h^{(k)}_{n}\}_{n=1}^{N_\omega}$ can be found. As a result, the distribution of $h$ can be explicitly described. Here XGBoost regressor is used as our simulator due to its great prediction accuracy \cite{chen2016xgboost}.

\begin{algorithm}
	\begin{small}
		{
			\caption{Data augmentation}
			\label{tab_algorithm1}
			\begin{tabular}{p{0.35cm}p{7.6cm}}
				$01$&\textbf{Simulator training:} train a regressor as our simulator based on the historical dataset $\{\bm x_n, \bm \omega_n, h_n\}_{n \in \mathcal{N}}$. Its input and output are $(\bm x, \bm \omega)$ and $h(\bm x, \bm \omega)$, respectively; \\ 
				$02$& \textbf{For} $k \in \mathcal{K}= [1,2,\cdots,K]$ \\
				$03$&\begin{adjustwidth}{3mm}{0cm}Randomly select one $\bm x$ and multiple $\bm \omega$ from the historical dataset to construct different pairs, and record them as $\{(\bm x^{(k)}, \bm \omega_n^{(k)})\}_{n=1}^{N_{\omega}}$ ($N_\omega$ is the number of $\bm \omega$ we chosen); \end{adjustwidth}\\
				$04$&\begin{adjustwidth}{3mm}{0cm} Give the above pairs as inputs to the simulator to predict $h$, and record the predictions as $\{\hat h^{(k)}_n\}_{n=1}^{N_{\omega}}$; \end{adjustwidth}\\
				$05$ &\textbf{End for}\\
				$06$ &\textbf{Data collection:} Collect the generated pairs and predictions to construct a new dataset, i.e., $\{\bm x^{(k)}, \bm \omega_n^{(k)},\hat h^{(k)}_n \}_{n=1}^{N_{\omega}}, \forall k \in \mathcal{K}$. By removing $\bm \omega_n^{(k)}$, we can get the ideal training set of the quantile-MLP, i.e., $\{\bm x^{(k)},\hat h^{(k)}_n \}_{n=1}^{N_{\omega}}, \forall k \in \mathcal{K}$.\\
		\end{tabular}} 
	\end{small}
\end{algorithm}

\subsubsection{Calibration}
According to (\ref{eqn_quantile_JCC})-(\ref{eqn_quantile_JCC2}), if the quantile-MLP overestimates the target quantile, i.e., $\mathcal{\hat Q}^{1-\epsilon} (\bm x) \geq \mathcal{Q}_{\bm \omega}^{1-\epsilon}(h(\bm x, \bm \omega))$, then only the optimality of solutions will be harmed but the feasibility can be maintained. However, if the quantile-MLP underestimates the target quantile, i.e., $\mathcal{\hat Q}^{1-\epsilon} (\bm x) \leq \mathcal{Q}_{\bm \omega}^{1-\epsilon}(h(\bm x, \bm \omega))$, then the feasibility may not be guaranteed. Thus, underestimation is more harmful and should be avoided. Based on this observation, a calibration step is developed. In this step, we design a special positive calibration parameter $\rho$ to calibrate the prediction of the quantile-MLP $\mathcal{\hat Q}^{1-\epsilon}(\bm x)$ so that the harmful underestimation can be avoided. Specifically, parameter $\rho$ is designed as the maximum underestimation of the simulator in the historical dataset $\{\bm x_n, \bm \omega_n, h_n\}_{n \in \mathcal{N}}$:
\begin{align}
\rho = \max_{n \in \mathcal{N}} \{h_n - \hat h_n\}, \label{eqn_rho_define}
\end{align}
where $\hat h_n$ is the prediction of the simulator based on $(\bm x_n, \bm \omega_n)$.
According to (\ref{eqn_rho_define}), we must have:
\begin{align}
\hat h_n  + \rho \geq h_n, \quad \forall n \in \mathcal{N}.
\end{align}
Then, in most cases, we have:
\begin{align}
&\mathcal{Q}_{\bm \omega}^{1-\epsilon}(\hat h + \rho) = \mathcal{Q}_{\bm \omega}^{1-\epsilon}(\hat h) + \rho \geq \mathcal{Q}_{\bm \omega}^{1-\epsilon}(h(\bm x, \bm \omega)), \label{eqn_overestimation}
\end{align}
where ``$=$" holds due to the translation invariance of the quantile. By regrading $\hat h$ as the training label to train the quantile-MLP, Eq. (\ref{eqn_overestimation}) can be replaced by:
\begin{align}
\mathcal{\hat Q}^{1-\epsilon}(\bm x) + \rho \geq \mathcal{Q}_{\bm \omega}^{1-\epsilon}(h(\bm x, \bm \omega)), \label{eqn_overestimation_2}
\end{align}
where $\mathcal{\hat Q}^{1-\epsilon}(\bm x)$ is the prediction of the quantile $\mathcal{Q}_{\bm \omega}^{1-\epsilon}(\hat h)$ given by the quantile-MLP. Note this prediction can achieve desirable accuracy because we can generate sufficient samples of $\hat h$ based on the simulator built in the data augmentation step. 
Finally, the quantile-based form of the joint chance constraint, (\ref{eqn_quantile_JCC}) can be innerly approximated by:
\begin{align}
 \mathcal{\hat Q}^{1-\epsilon}(\bm x) + \rho \leq 0.
\label{eqn_quantile_constraint3}
\end{align}
According to (\ref{eqn_overestimation_2}), Eq. (\ref{eqn_quantile_constraint3}) is an inner approximation of (\ref{eqn_quantile_JCC}). Thus, the aforementioned harmful underestimation can be avoided and the feasibility of solutions can be guaranteed.

\subsection{Whole procedure of the proposed model}
By applying the two pre-processing steps, \textbf{P2} can be replaced by the following surrogate model \textbf{P3}:
\begin{align} 
	&\min_{\bm \lambda, G} \quad \mathbb E_{\bm \omega} (G) \tag{$\textbf{P3}$},\\
	&\begin{array}{r@{\quad}r@{}l@{\quad}l}
		\text{s.t.:} &&\text{Eqs. (\ref{eqn_x_define}), (\ref{eqn_input}), (\ref{eqn_output}), (\ref{eqn_balance2}), (\ref{eqn_input_loss}), (\ref{eqn_output_loss})-(\ref{eqn_reformulation_loss}), and (\ref{eqn_quantile_constraint3}).} 
	\end{array} \notag
\end{align}

\begin{figure}
	\centering
	 			\vspace{-4mm}
	\includegraphics[width=1\columnwidth]{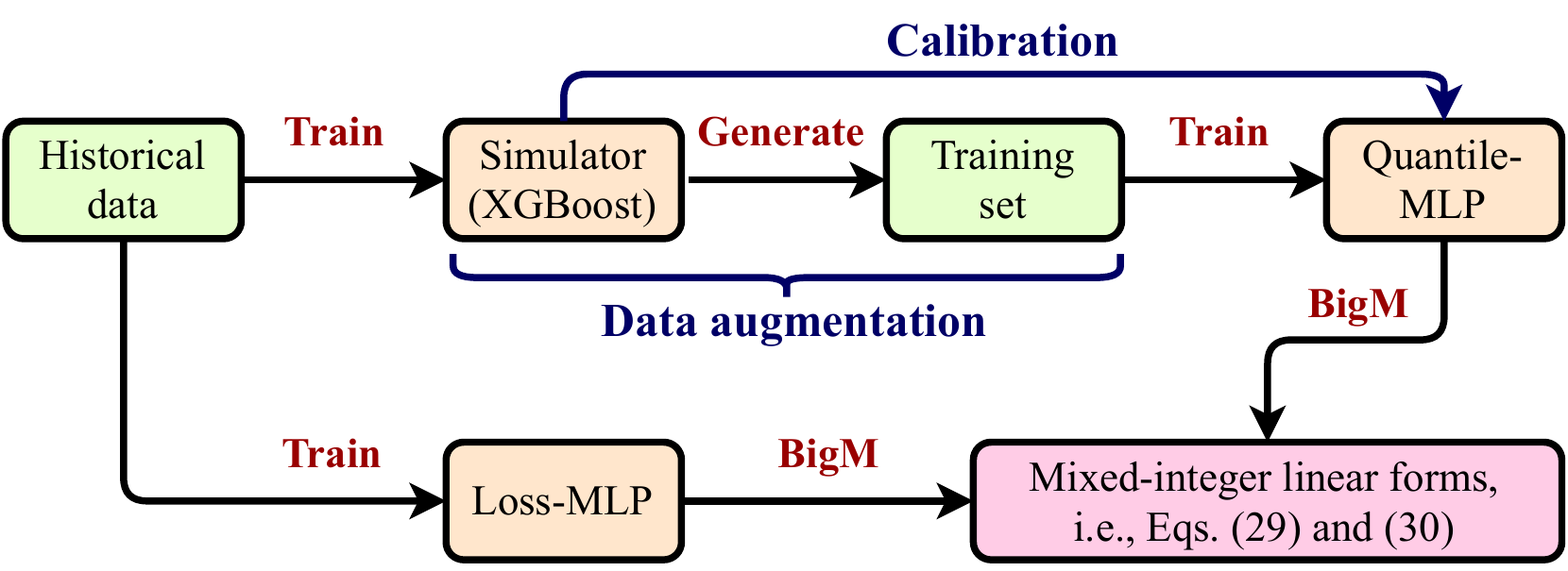}\vspace{-2mm}
	\caption{The whole procedure to establish the proposed surrogate model.
	}
	\label{fig_procedure}
	 		\vspace{-2mm}
\end{figure}

Fig. \ref{fig_procedure} illustrates the whole procedure for establishing the proposed learning-based surrogate model. Specifically, we first leverage historical data to train a simulator for data augmentation. Then, the quantile MLP is trained to predict the target quantile $\mathcal{Q}_{\bm \omega}^{1-\epsilon}(h(\bm x, \bm \omega))$. A calibration step is further applied to the quantile-MLP for improving the feasibility of solutions. Meanwhile, another MLP, i.e., the loss-MLP, is trained based on the historical data to predict the expected power loss. Finally, by reformulating these two MLPs into solvable mixed-integer forms, the proposed surrogate model can be established to replicate the JCC-OPF problem without the network parameters.

In the proposed model, three regressors are trained. For convenience, we summarize their inputs, outputs, and training sets in Table \ref{tab_regressor}. The required historical samples include nominal active/reactive power injections, bus voltages, branch currents, and nominal/actual available DGs' outputs in the past.

\begin{table}
\small
\centering
\caption{Descriptions of the trained regressors}
\vspace{-2mm}
\begin{tabular}{cccc}
\hline
\textbf{Regressors} & \textbf{Inputs} & \textbf{Outputs} & \textbf{Training set} \\ \hline
Simulator           & $(\bm x,\bm \omega)$           & $\hat{h}$           & Historical dataset       \\
Quantile-MLP            & $\bm x$               &      ${\hat{\mathcal{Q}}}^{1-\epsilon}(\bm x)$     & Augmented dataset        \\
Loss-MLP             & $\bm x$               & $\hat p^\text{loss}$        & Historical dataset       \\ \hline
\end{tabular} \label{tab_regressor}
\vspace{-4mm}
\end{table}

\section{Case study} \label{sec_case}
\subsection{Simulation set up}
\begin{figure}
		\vspace{-4mm}
	\subfigbottomskip=-4pt
	\subfigcapskip=-4pt
	\centering
	\subfigure[]{\includegraphics[width=0.85\columnwidth]{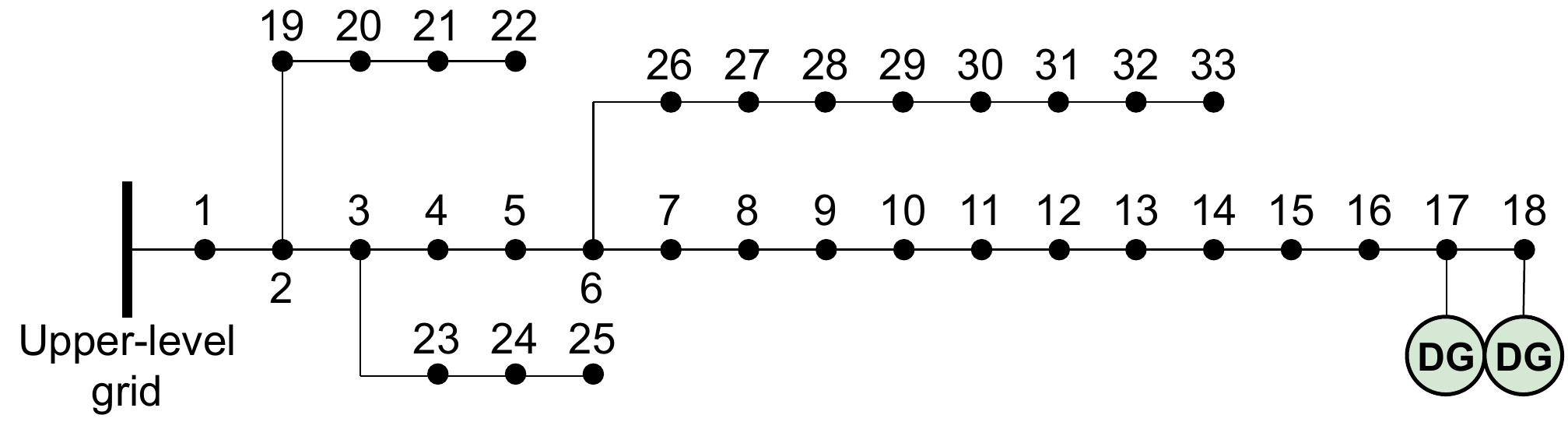}}
	\subfigure[]{\includegraphics[width=1\columnwidth]{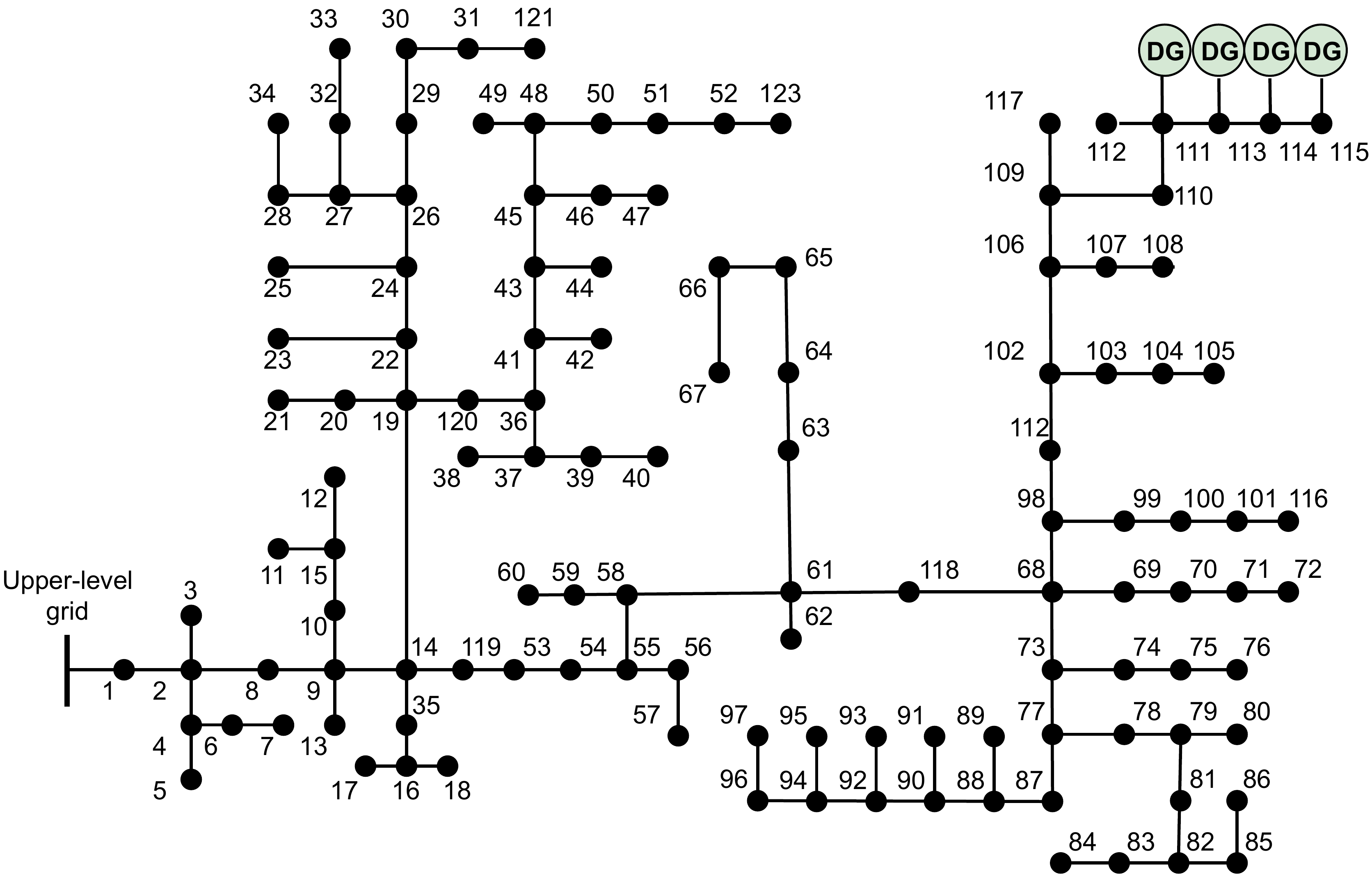}}
	\vspace{-4mm}
 	\caption{Structures of the (a) 33-bus system and (b) 123-bus system.}
	\label{fig_33bus}
	\vspace{-4mm}
\end{figure}

We implement two different case studies based on on the IEEE 33- and 123-bus systems, respectively. The first case study has two DGs, while the second one contains four DGs. There structures are illustrated in Fig. \ref{fig_33bus}. The slack bus voltage, i.e., $V_1$, in the two case studies are 12.66kV and 4.16kV, respectively. Other parameters used in the simulations are summarized in Table \ref{tab_parameter}.

\begin{table}
	\small
	\centering
	\caption{Parameters in simulations}
	\vspace{-2mm}
	\begin{threeparttable} 
	\begin{tabular}{ccccc}
		\hline
		\rule{0pt}{11pt}		
		Case studies & Parameters & Value &Parameters & Value\\
		\hline
		\rule{0pt}{10pt}
		\multirow{3}{*}{33-bus} & $\overline{G}_i^\text{DG}$ &  2MW & $V_\text{i, max}$& 1.1 p.u.\\
		&$I_\text{b, max}$ &  0.249kA & $V_\text{i, min}$& 0.9 p.u. \\
		& $\bm \phi$ & 0 & &\\
		\hline
		\multirow{2}{*}{123-bus} & $\overline{G}_i^\text{DG}$ &  1.5MW & $V_\text{i, max}$& 1.1 p.u.\\
		&$I_\text{b, max}$ &  0.65kA & $V_\text{i, min}$& 0.9 p.u. \\
		& $\bm \phi$ & 0.33 & &\\
		\hline
	\end{tabular}\label{tab_parameter}
	\end{threeparttable}
	\vspace{-2mm}
\end{table}

We conduct power flow simulations based on Pandapower, a power system simulation toolbox in Python environment \cite{8344496}, to generate the historical data. In Pandapower, the power flow calculation is based on the full AC power flow model. During the simulations, we first randomly generate 10,000 pairs of nominal bus power injections and uncertain levels of DGs' outputs, i.e., $(\bm x, \bm \omega)$. Here the nominal power injection $\bm x$ is generated by a uniform distribution between its minimum and maximum allowable values. Based on these pairs, we can calculate the actual power injections on each bus, and then the bus voltages $\bm V$ and branch currents $\bm I$ can be calculated by Pandapower. With $\bm V$ and $\bm I$, the power loss $p^\text{loss}$ and maximum constraint violation $h(\bm x, \bm \omega)$ can be obtained based on (\ref{eqn_loss})-(\ref{eqn_h_define}). Then, following \textbf{Algorithm} \ref{tab_algorithm1}, we conduct the data augmentation to generate the training set for the quantile-MLP. The two parameters $N_{\omega}$ and $K$ in \textbf{Algorithm} \ref{tab_algorithm1} are set as 1000 and 100, respectively.

To demonstrate the generalization ability of the proposed model, different distributions are used to generate the samples of the uncertain level $\bm \omega$, as follows:
\begin{enumerate}
	\item \textbf{Case 1}: The samples of $\bm \omega$ are generated by a Gaussian distribution, i.e., $\bm \omega \sim \text{Gaussian}(0,0.1)$;
	\item \textbf{Case 2}: The samples of $\bm \omega$ are generated based on a Beta distributed uncertainty $\bm \omega^{'}$, i.e., $\bm \omega = \kappa^\text{Beta}(\bm \omega^{'} - \bm \mu^\text{beta})$, where $ \bm \omega^{'} \sim \text{Beta}(2,6)$;
	\item \textbf{Case 3}: The samples of $\bm \omega$ are generated based on a Weibull distributed uncertainty $\bm \omega^{'}$, i.e., $\bm \omega = \kappa^\text{Weibull}(\bm \omega^{'} - \bm \mu^\text{Weibull})$, where $\bm \omega^{'} \sim \text{Weibull}(1,5)$.
\end{enumerate}
The scaling factor $\kappa^\text{Beta}/\kappa^\text{Weibull}$ is designed to make the magnitudes of the generated $\bm \omega$ more realistic, while the constant $\bm \mu$ is set as the expectation of $\bm \omega^{'}$ to make the expectation of $\bm \omega$ keep at zero. All these samples have been uploaded to \cite{samples2022}. 


All numerical experiments are implemented on an Intel(R) 8700 3.20GHz CPU with 16 GB memory. The quantile-MLP and loss-MLP are established based on Pytorch. Dropout is also applied to mitigate the overfitting \cite{srivastava2014dropout}. Problem \textbf{P3} is built by CVXPY and solved by GUROBI.

\subsection{Benchmarks}
To demonstrate the superiority of the proposed model, we introduce the following three benchmarks:
\begin{enumerate}
	\item \textbf{B1}: Linearized DistFlow model used in \cite{8515118,8688432,8528856} combined with the scenario approach;
	\item \textbf{B2}: SOCP relaxation of AC OPF model used in \cite{8060613, 8626040} combined with the scenario approach;
	\item \textbf{B3}: Risk-neutral full AC OPF model.   
\end{enumerate}
In \textbf{B1} and \textbf{B2}, the scenario approach used in \cite{8060613, 8626040} is employed to handle the intractable joint chance constraint (\ref{eqn_JCC}). \textbf{B3} is a non-convex problem, which is directly solved by Pandapower.
Note that \textbf{B1}-\textbf{B3} are based on power flow models, in which the network parameters are assumed known. 

\subsection{Case study based on the 33-bus system} \label{sec_33bus}
To verify the effectiveness of the proposed model, we compare the average utilization rates of DG, maximum violation probabilities of the joint chance constraint  (\ref{eqn_JCC}), energy purchasing from the upper-level grid, and solving times of different models based on the 33-bus test system.
Both the maximum violation probabilities and energy purchasing are obtained by Monte Carlo Simulations based on Pandpower to make the comparison fair\footnote{We first solve each method to the solution. Then, by giving these solutions and uncertainty samples as inputs to Pandapower, the actual maximum violation probabilities and energy purchasing can be calculated.}. The neuron numbers of the quantile-MLP are set as (25, 25, 25), i.e., three hidden layers with 25 neurons in each layer, while the neuron numbers of the loss-MLP are set as (10, 10, 10). 

\subsubsection{Case 1}

\begin{figure}
		\vspace{-4mm}
	\subfigbottomskip=-4pt
	\subfigcapskip=-4pt
	\centering
	\subfigure[]{\includegraphics[width=0.49\columnwidth]{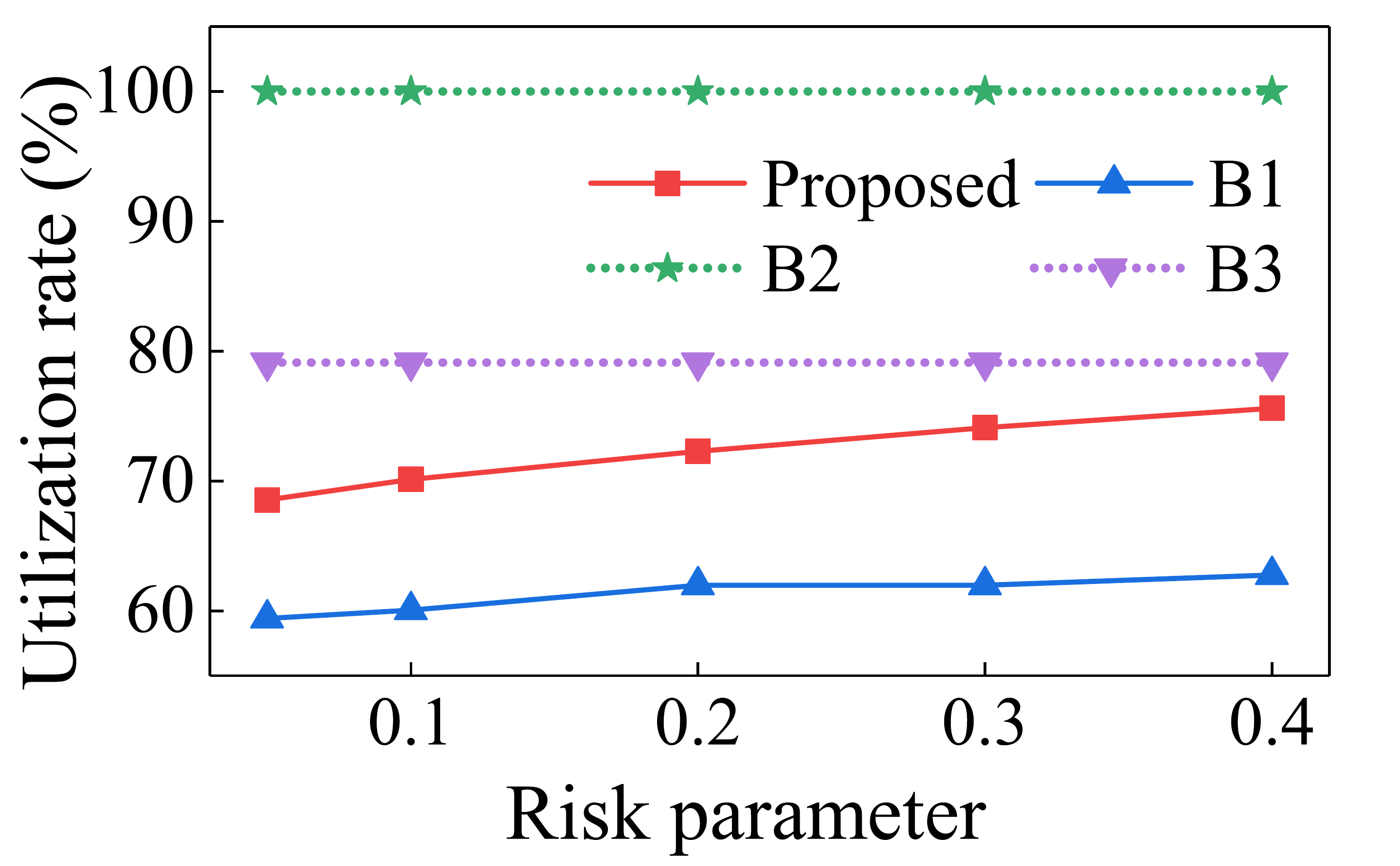}}
	\subfigure[]{\includegraphics[width=0.49\columnwidth]{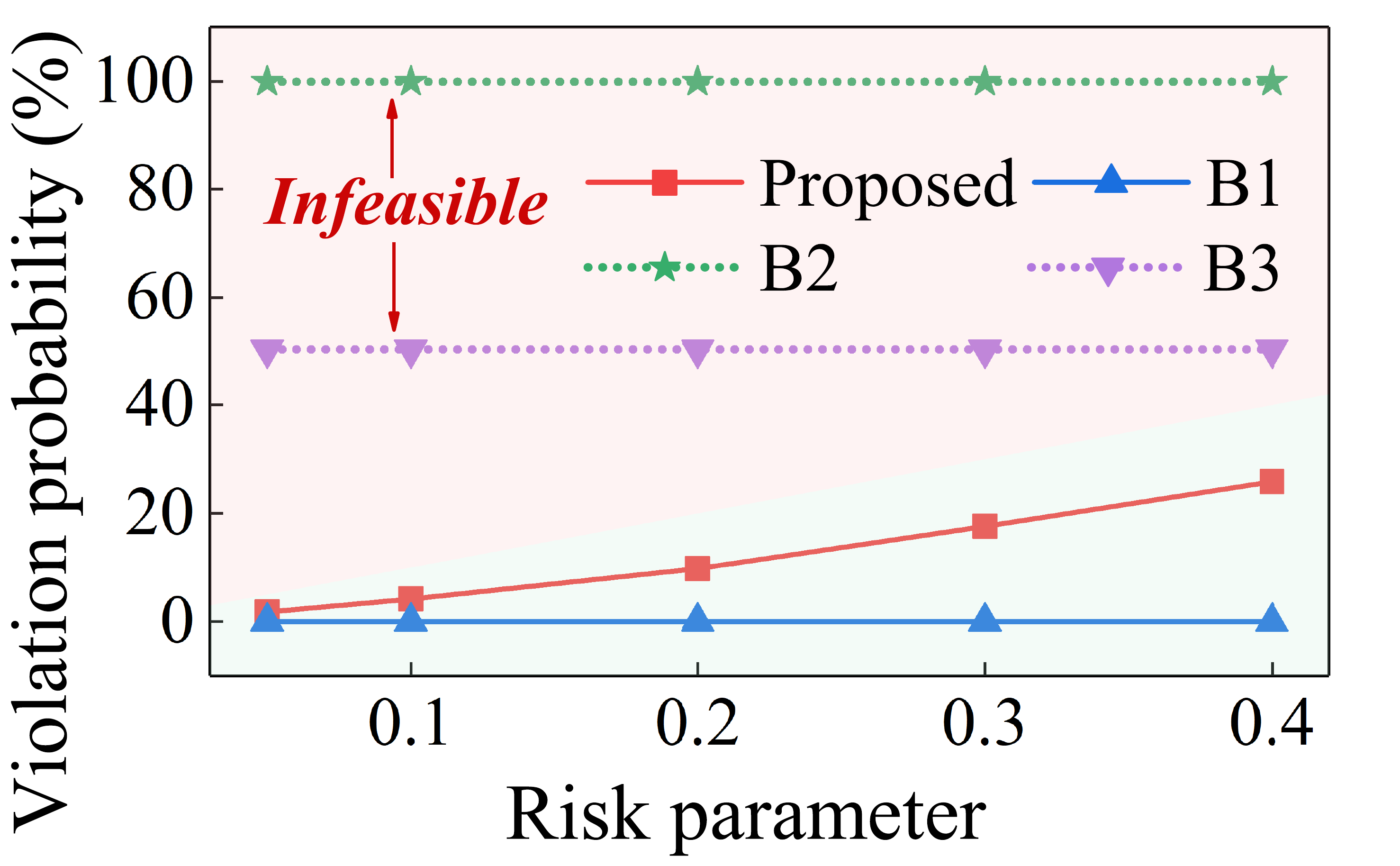}}
	\subfigure[]{\includegraphics[width=0.49\columnwidth]{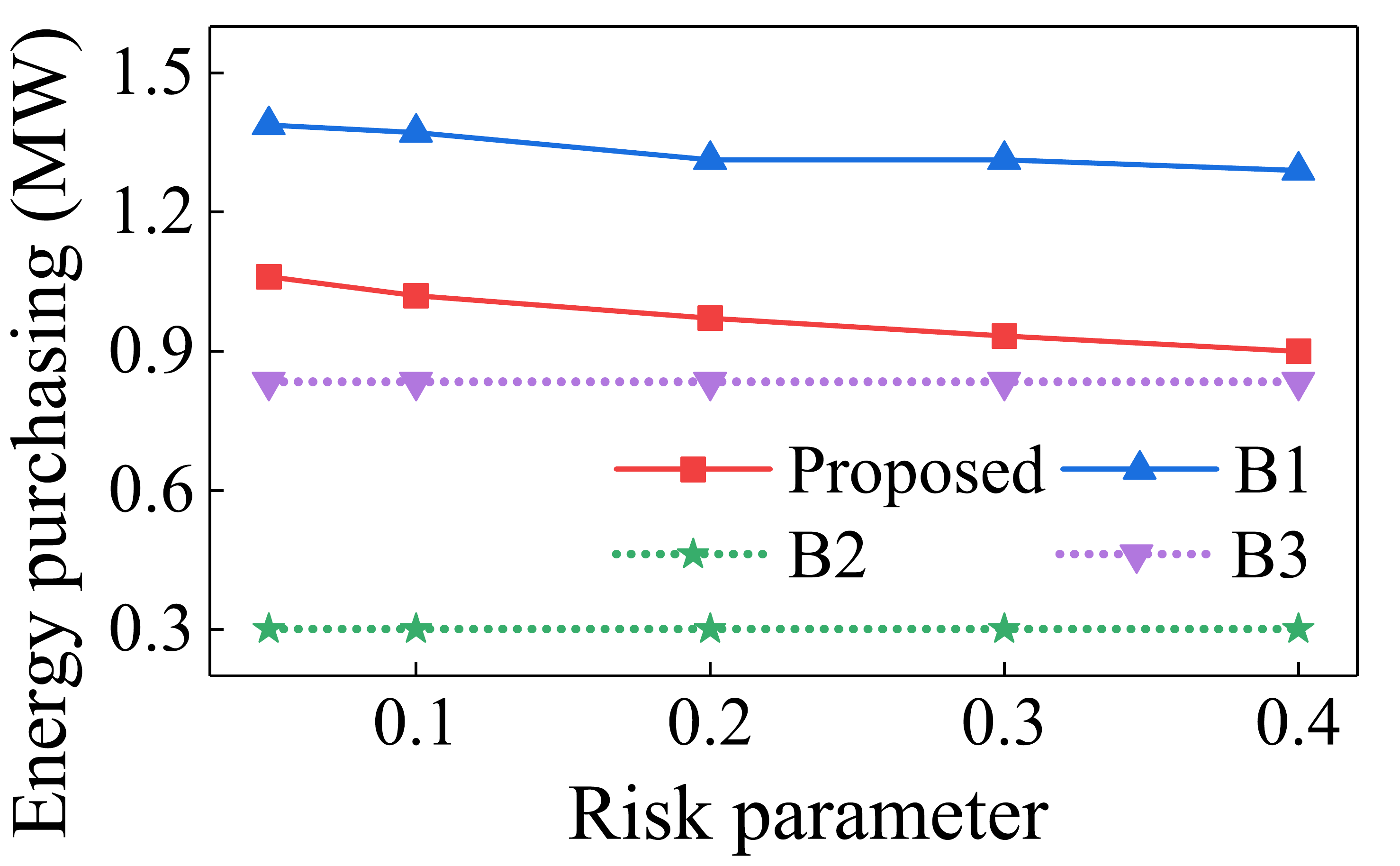}}
	\subfigure[]{\includegraphics[width=0.49\columnwidth]{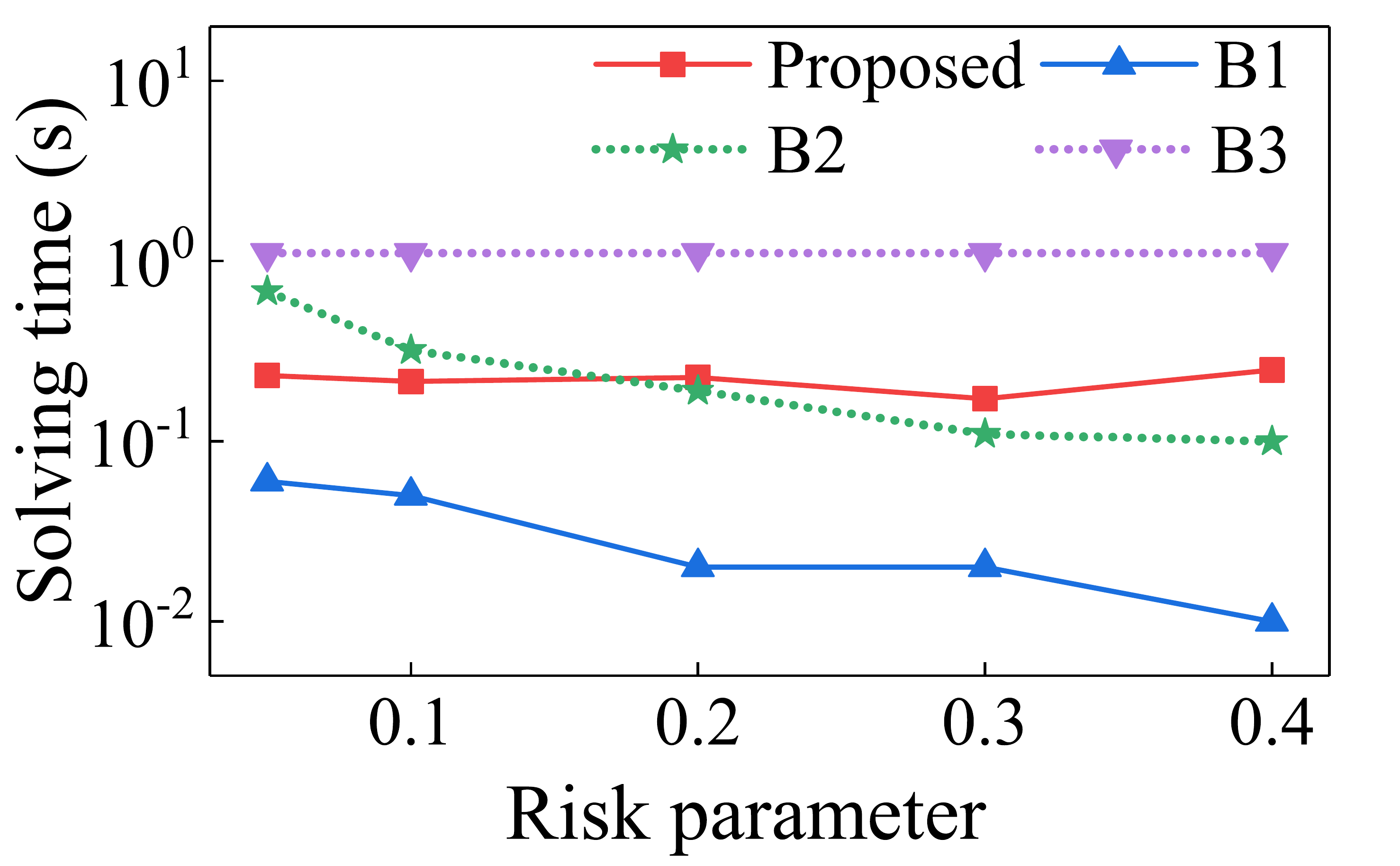}}
	\vspace{-4mm}
 	\caption{Results of (a) average utilization rates of DG, (b) maximum violation probabilities of the joint chance constraint  (\ref{eqn_JCC}), (c) energy purchasing from the upper-level grid, and (d) solving times in Case 1 (the uncertainty $\bm \omega$ follows Gaussian distribution). Dot lines represent infeasible results. In (b), the green and red areas denote the feasible and infeasible regions, respectively.}
	\label{fig_results_Gaussian}
	\vspace{-4mm}
\end{figure}
Figure \ref{fig_results_Gaussian} compares the results of different models in \textbf{Case 1}.
Among all models, the linearized DistFlow model, \textbf{B1}, derives the most conservative results, while its maximum violation probability and utilization rate of DG are the lowest. Since \textbf{B1} ignores voltage drops on branches, it overestimates bus voltages \cite{19266}. Considering that promoting the integration of DG will increase bus voltages, less DG can be utilized in \textbf{B1} because it must ensure the overestimated bus voltages to be smaller than the corresponding upper bound. Nevertheless, its solution is always feasible for the joint chance constraint  (\ref{eqn_JCC}). Conversely, the SOCP relaxation \textbf{B2} shows very poor feasibility (its maximum violation probability approaches 100\%), although it achieves the highest utilization rate of DG and lowest energy purchasing. This is because reverse power flows occur in the system, which makes its SOCP relaxation inexact. The risk-neutral model \textbf{B3} also fails to meet the joint chance constraint (\ref{eqn_JCC}) since it directly ignores the impacts of uncertainties. Thus, both \textbf{B2} and \textbf{B3} are not applicable to distribution systems with high DG penetration. 
The proposed model can always ensure the feasibility of solutions, and its energy-efficiency is better than that of \textbf{B1}. Moreover, unlike the three benchmarks \textbf{B1}-\textbf{B3}, the proposed model only need historical data to train MLPs but does not require the network parameters to build power flow models. Although its computational performance is worse than \textbf{B1} due to the binary variables introduced by reformulating MLPs, the solving time is still acceptable (around 0.3s). These results demonstrate the great feasibility and optimality of the proposed model.


\subsubsection{Case 2}
Figure \ref{fig_results_Gaussian} compares the results of different models in \textbf{Case 2}, in which the uncertainty $\bm \omega$ follows Beta distribution. The results are very similar to those in \textbf{Case 1}: the energy-efficiency of \textbf{B1} is undesirable, while \textbf{B2} and \textbf{B3} can not guarantee the feasibility of solutions. In contrast, the proposed model can achieve desirable optimality and feasibility simultaneously with no need for the network parameters. The solving time of the proposed model still keeps around 0.3s, which is acceptable in practice.

\begin{figure}
		\vspace{-4mm}
	\subfigbottomskip=-4pt
	\subfigcapskip=-4pt
	\centering
	\subfigure[]{\includegraphics[width=0.49\columnwidth]{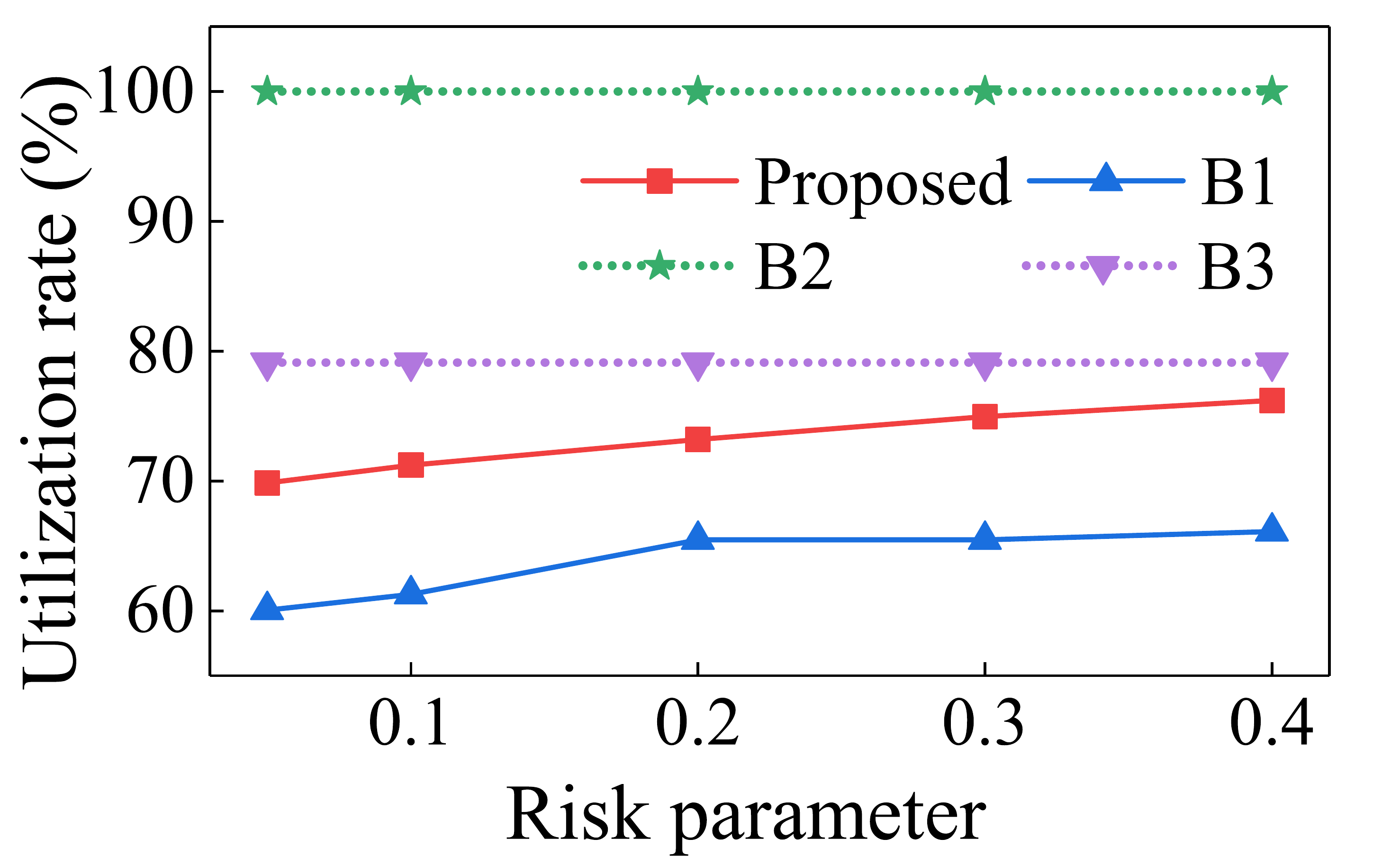}}
	\subfigure[]{\includegraphics[width=0.49\columnwidth]{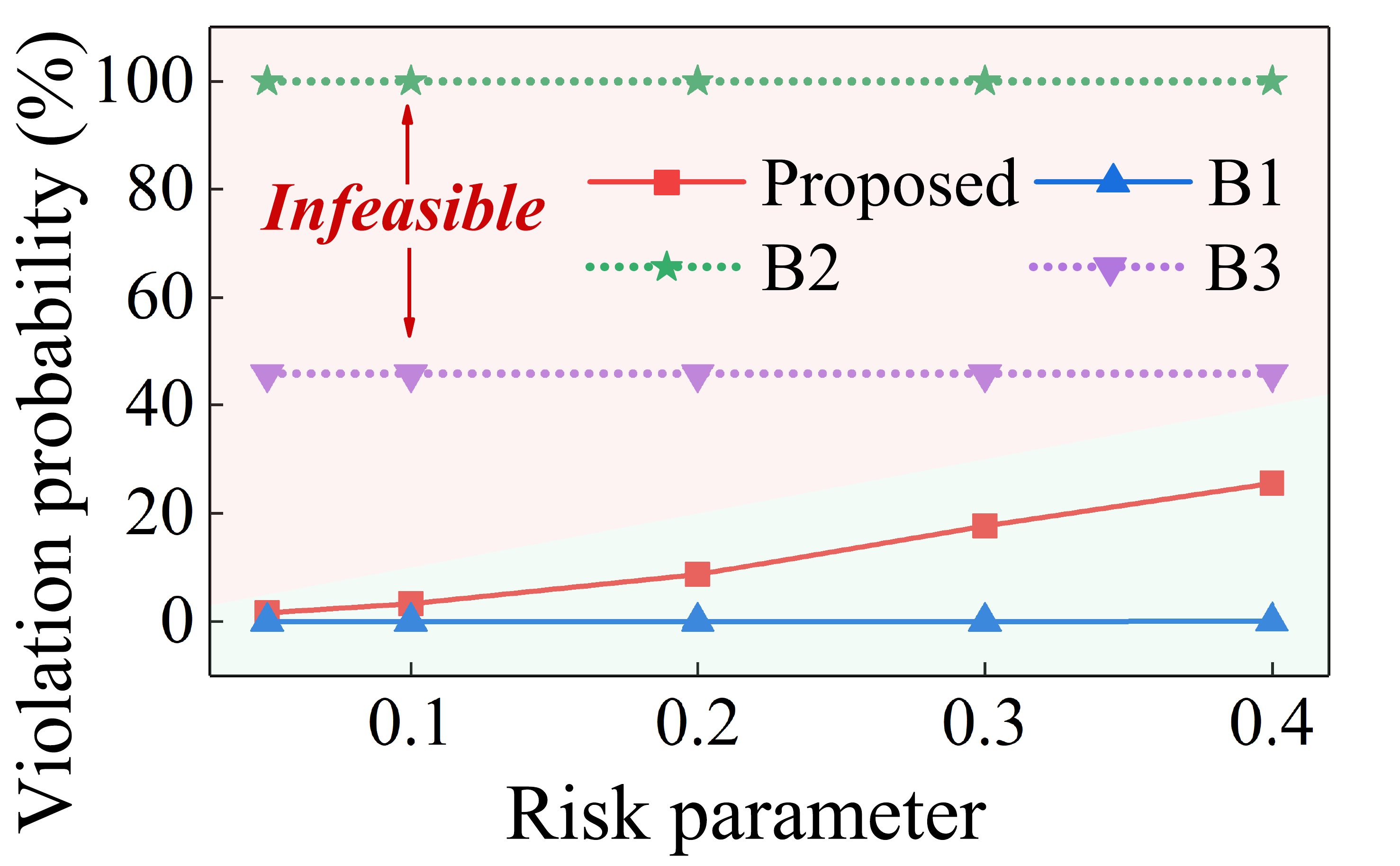}}
	\subfigure[]{\includegraphics[width=0.49\columnwidth]{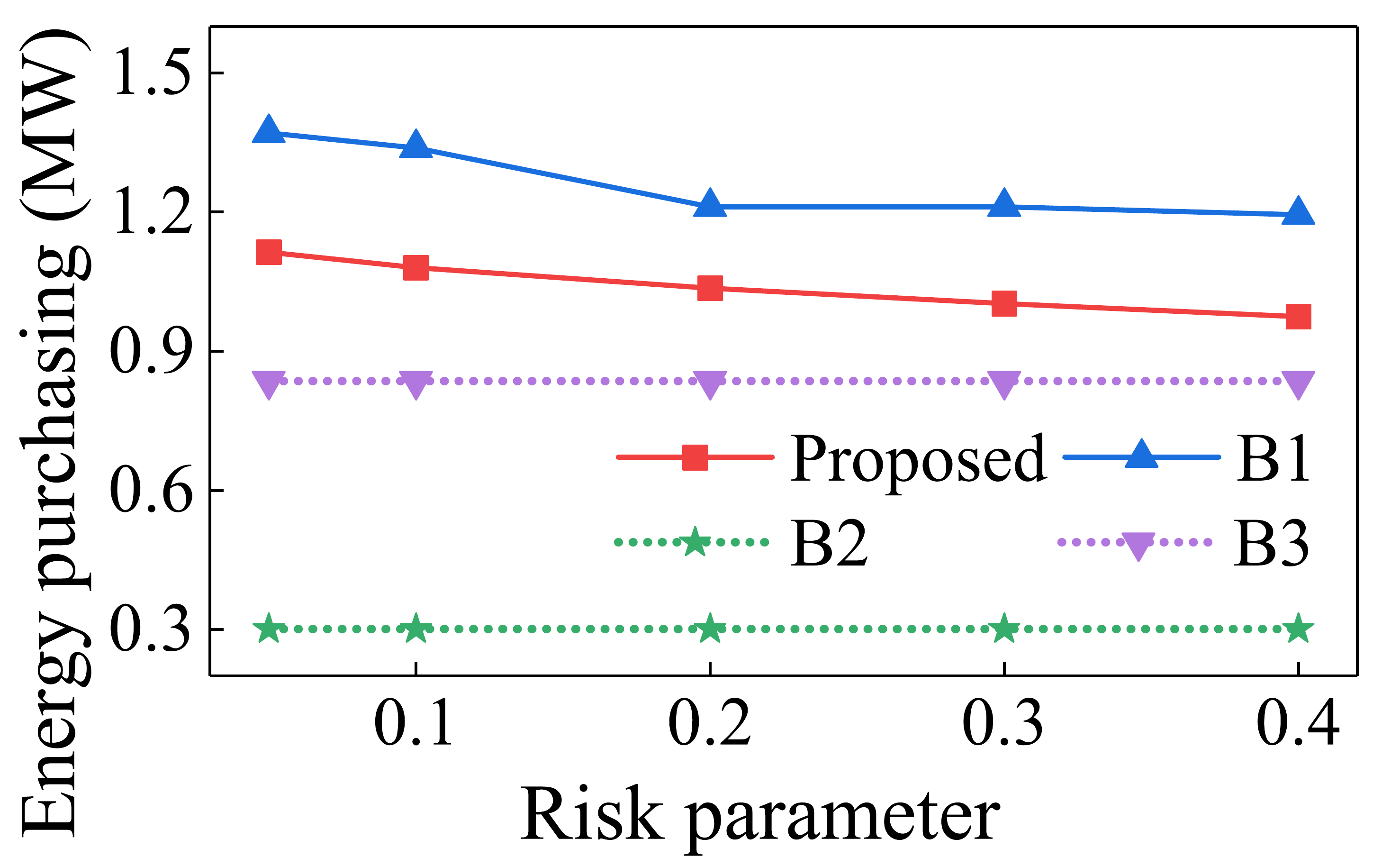}}
	\subfigure[]{\includegraphics[width=0.49\columnwidth]{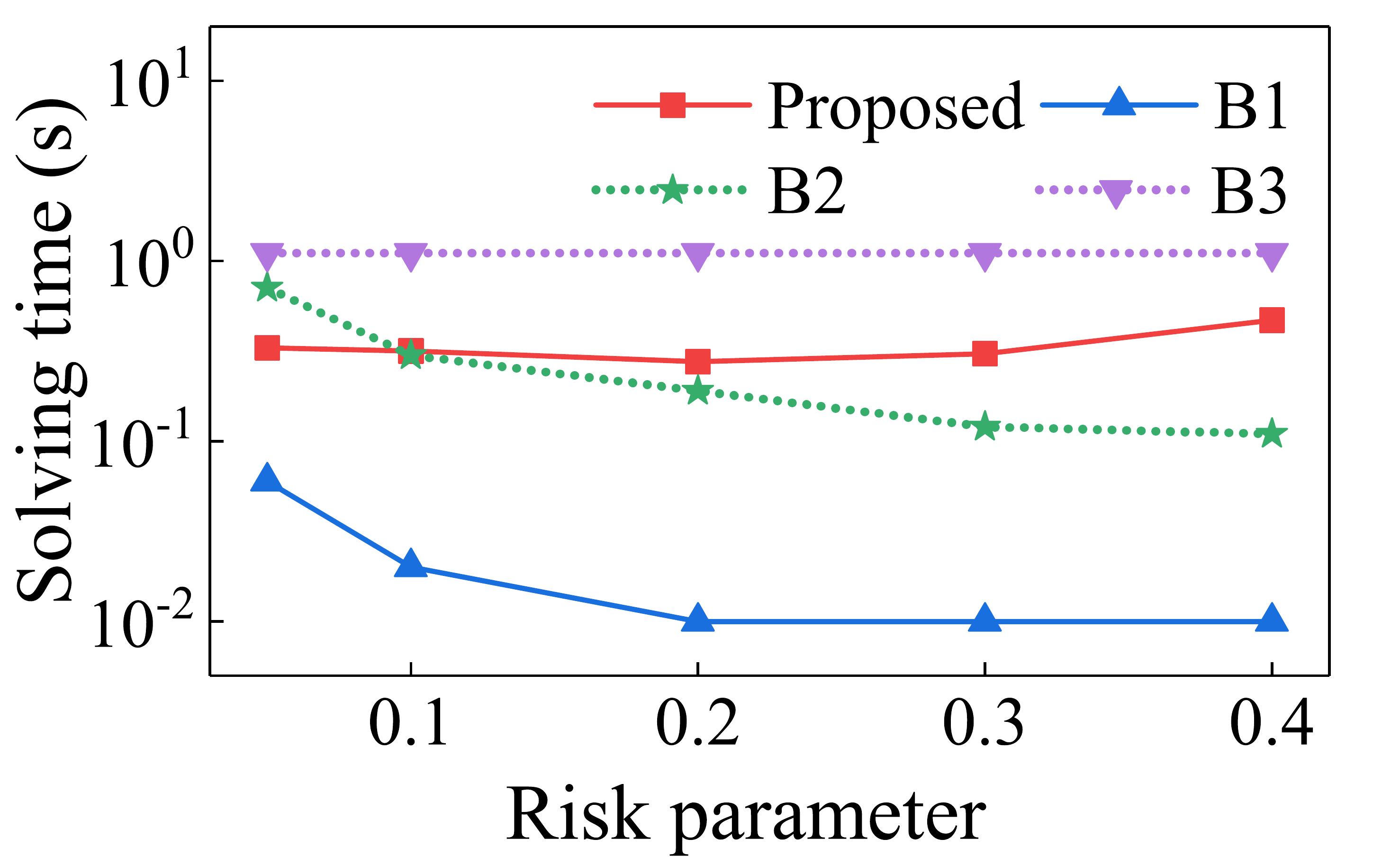}}
	\vspace{-4mm}
 	\caption{Results of (a) average utilization rates of DG, (b) maximum violation probabilities, (c) energy purchasing, and (d) solving times in Case 2 (the uncertainty $\bm \omega$ follows Beta distribution).}
	\label{fig_results_Beta}
	\vspace{-4mm}
\end{figure}

\subsubsection{Case 3}
Figure \ref{fig_results_Gaussian} illustrates the results of different models in \textbf{Case 3}, in which the uncertainty $\bm \omega$ follows Weibull distribution. Similarly, the proposed model achieves better energy-efficiency compares to \textbf{B1}. Moreover, it outperforms \textbf{B2} and \textbf{B3} on feasibility. 

In summary, the above three cases not only illustrate that the proposed model can achieve desirable optimality and feasibility without the topology but also demonstrate its excellent generalization performance for arbitrary uncertainties.

\begin{figure}
	\subfigbottomskip=-4pt
	\subfigcapskip=-4pt
	\centering
	\subfigure[]{\includegraphics[width=0.49\columnwidth]{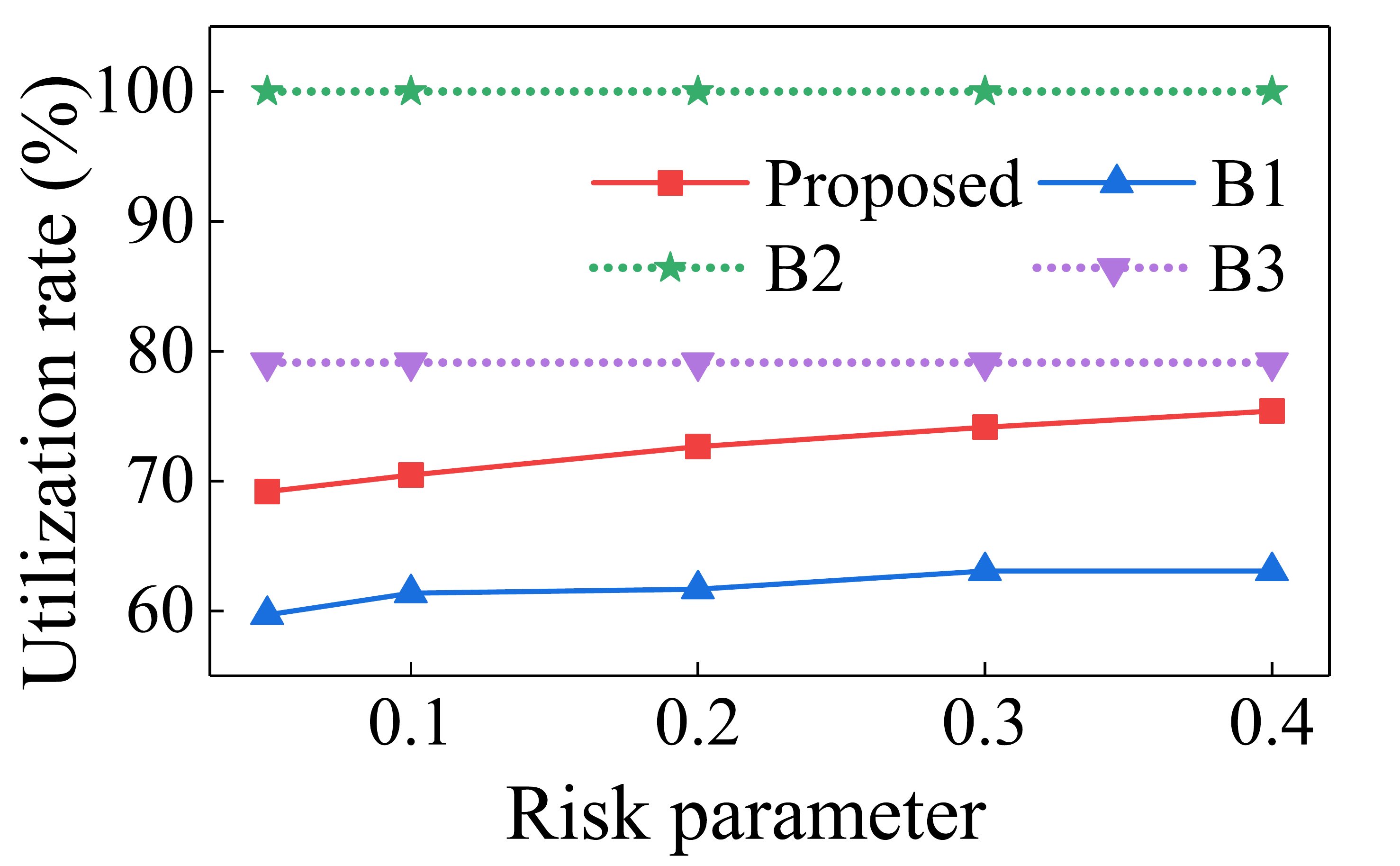}}
	\subfigure[]{\includegraphics[width=0.49\columnwidth]{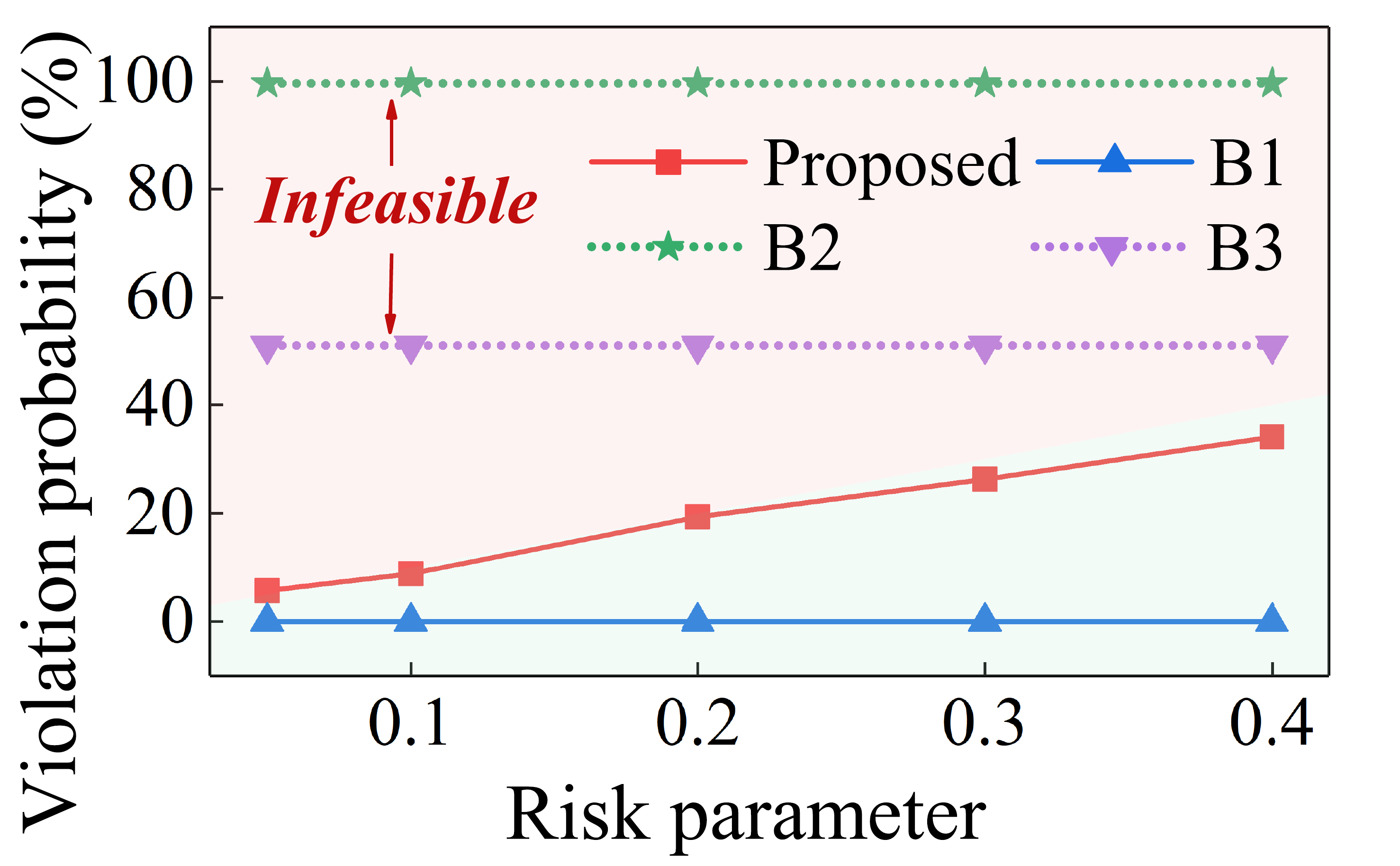}}
	\subfigure[]{\includegraphics[width=0.49\columnwidth]{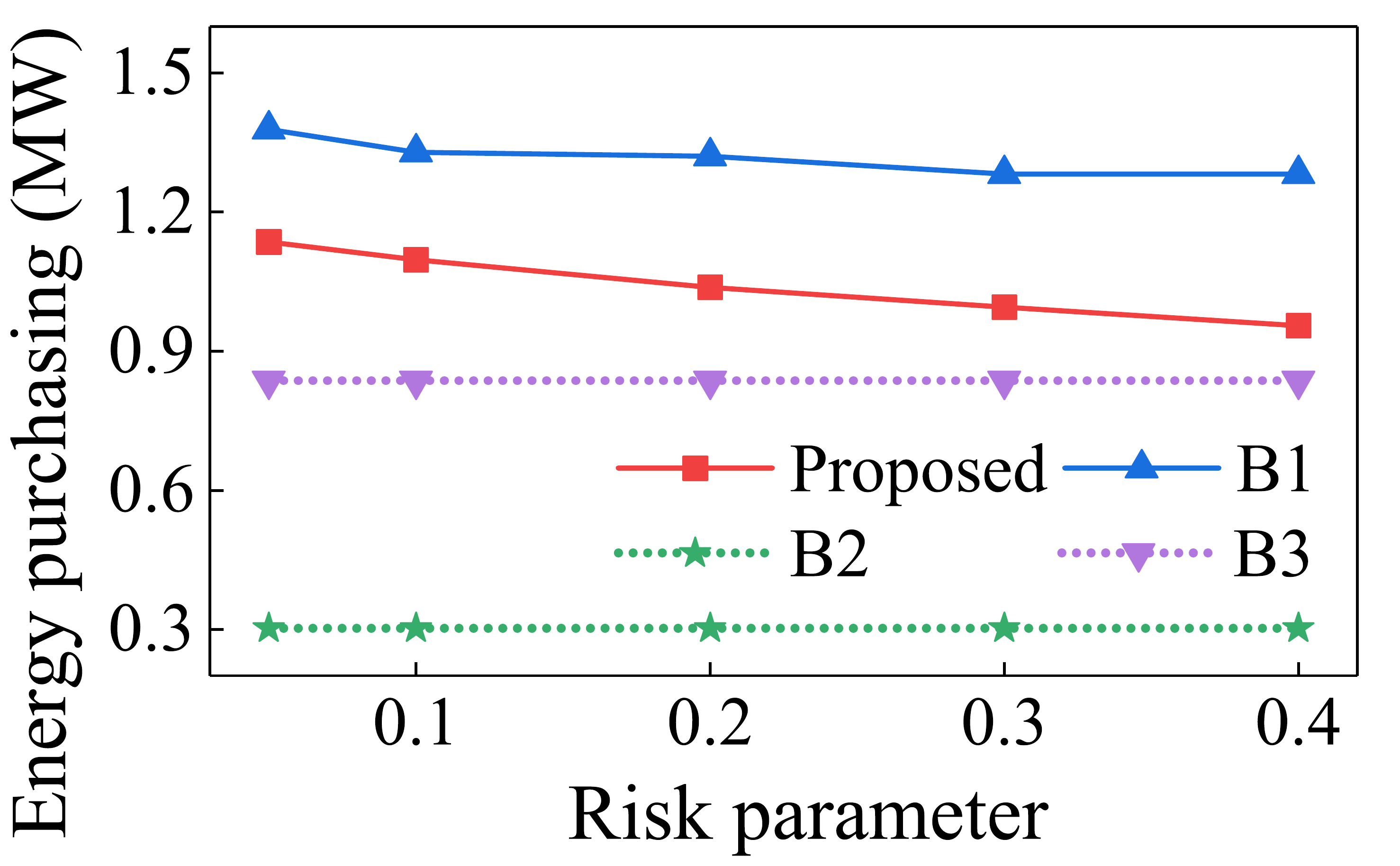}}
	\subfigure[]{\includegraphics[width=0.49\columnwidth]{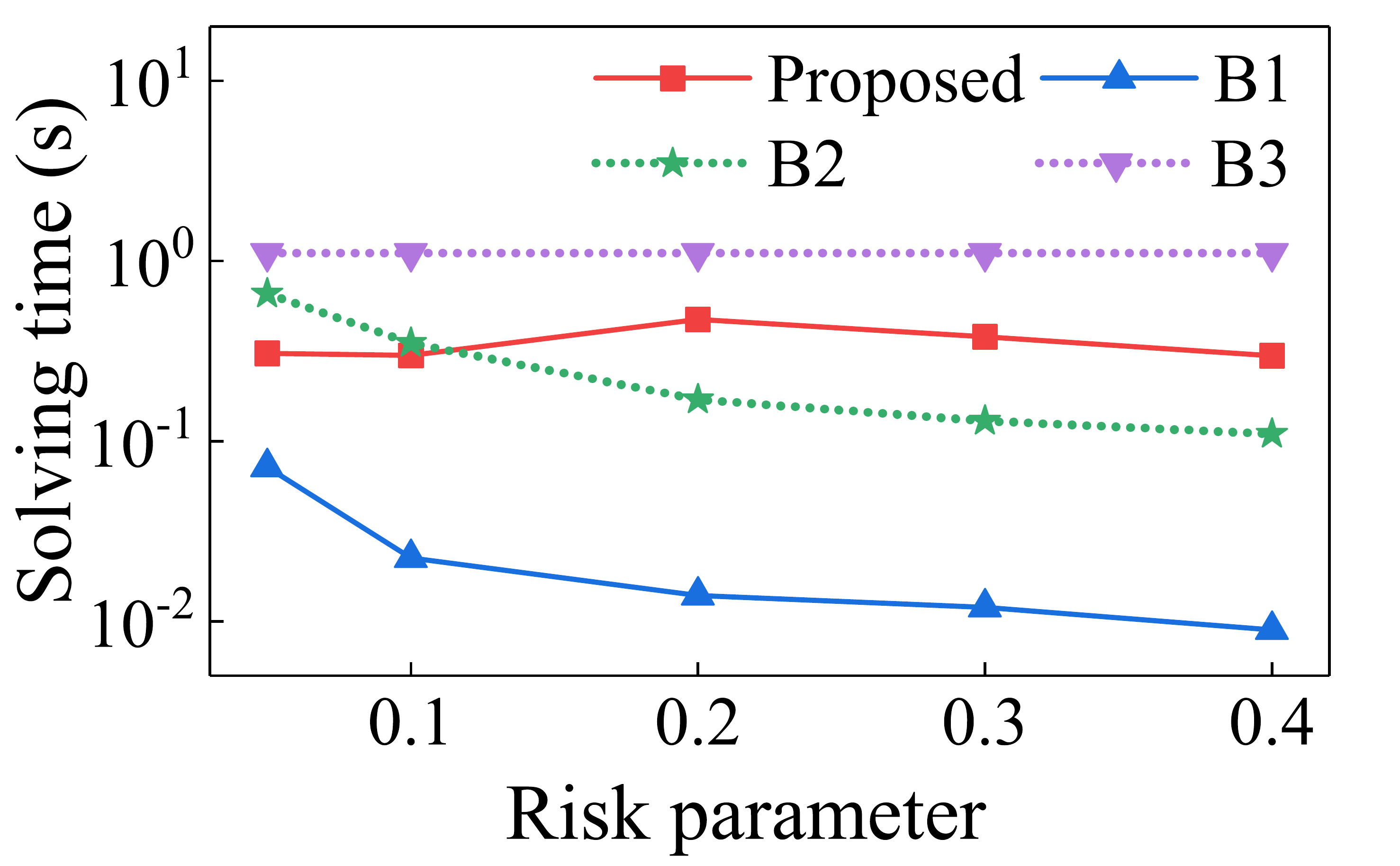}}
	\vspace{-4mm}
 	\caption{Results of (a) average utilization rates of DG, (b) maximum violation probabilities, (c) energy purchasing, and (d) solving times in Case 3 (the uncertainty $\bm \omega$ follows Weibull distribution). }
	\label{fig_results_Weibull}
	\vspace{-4mm}
\end{figure}

\subsection{Case study based on the 123-bus system}
We further conduct a case study based on the 123-bus system to better demonstrate the benefits of the proposed model. The neuron numbers of the quantile- and loss-MLPs are set as (30, 30, 30) and (10, 10, 10), respectively. The uncertainties are the same as those in \textbf{Case 3} (Weibull uncertainties). 

The results in Section \ref{sec_33bus} show that \textbf{B1} is overly conservative. However, this conservativeness may be contributed by either the power flow approximation (linearized DistFlow) or the JCC reformulation (scenario approach). To highlight that the linearized DistFlow model introduces unnecessary conservativeness, we modify the benchmark \textbf{B1} as \textbf{B1-SAA}, in which the intractable JCCs are handled by sample average approximation (SAA). SAA is a promising way to handle JCCs with excellent optimality, but it is also time-consuming because numerous binary variables have to be involved \cite{geng2019data}.

\begin{figure}
	\subfigbottomskip=-4pt
	\subfigcapskip=-4pt
	\centering
	\subfigure[]{\includegraphics[width=0.49\columnwidth]{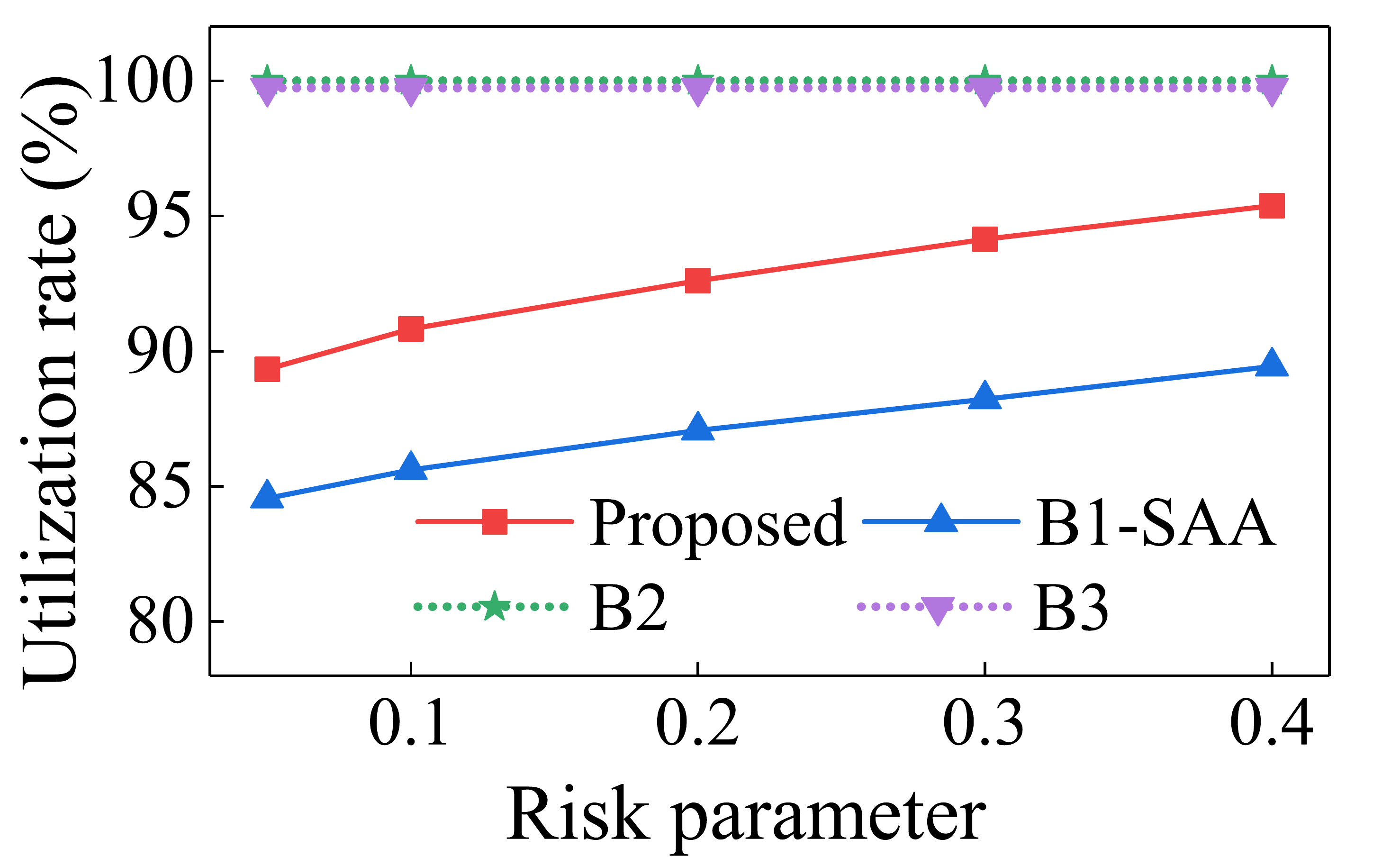}}
	\subfigure[]{\includegraphics[width=0.49\columnwidth]{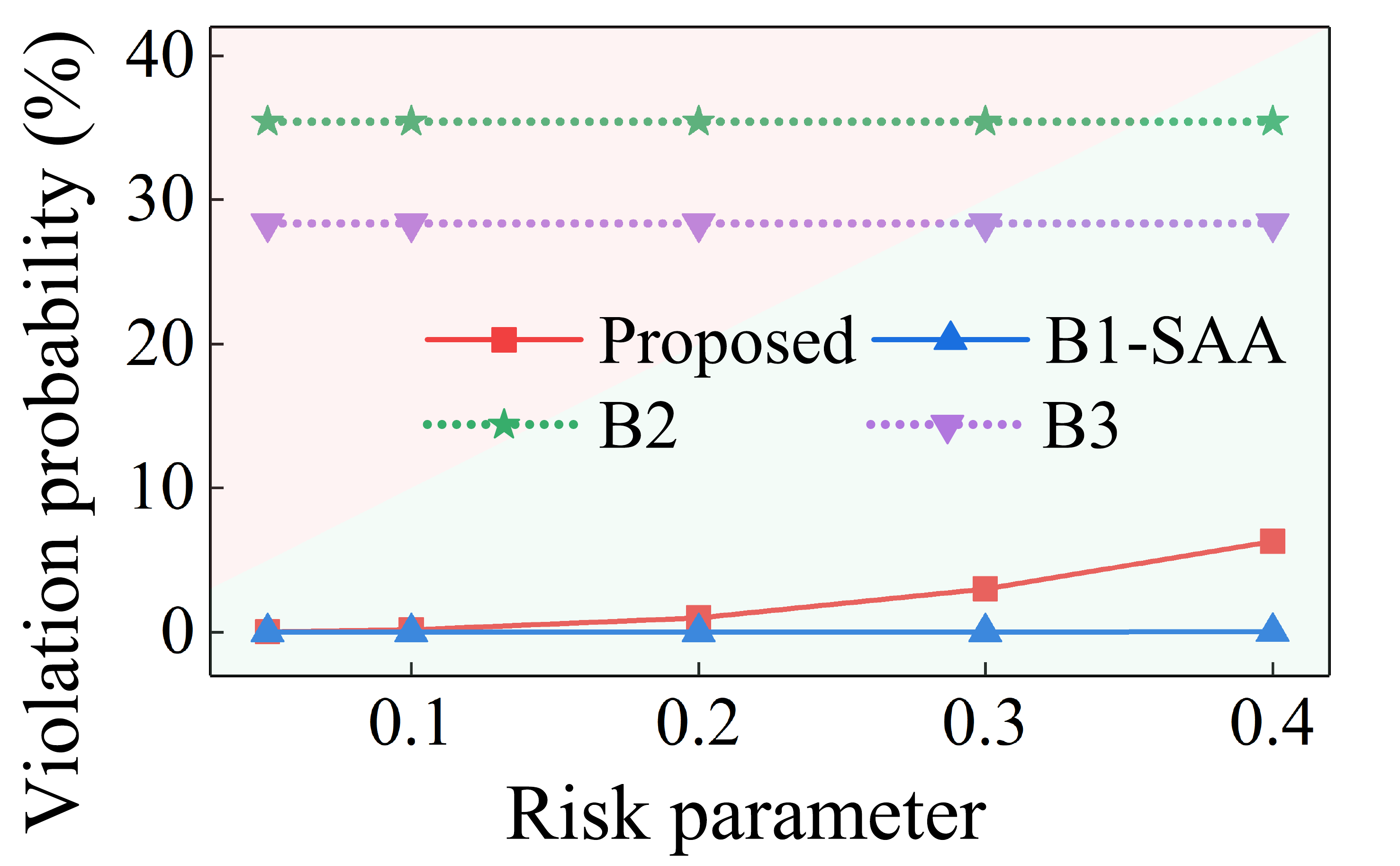}}
	\subfigure[]{\includegraphics[width=0.49\columnwidth]{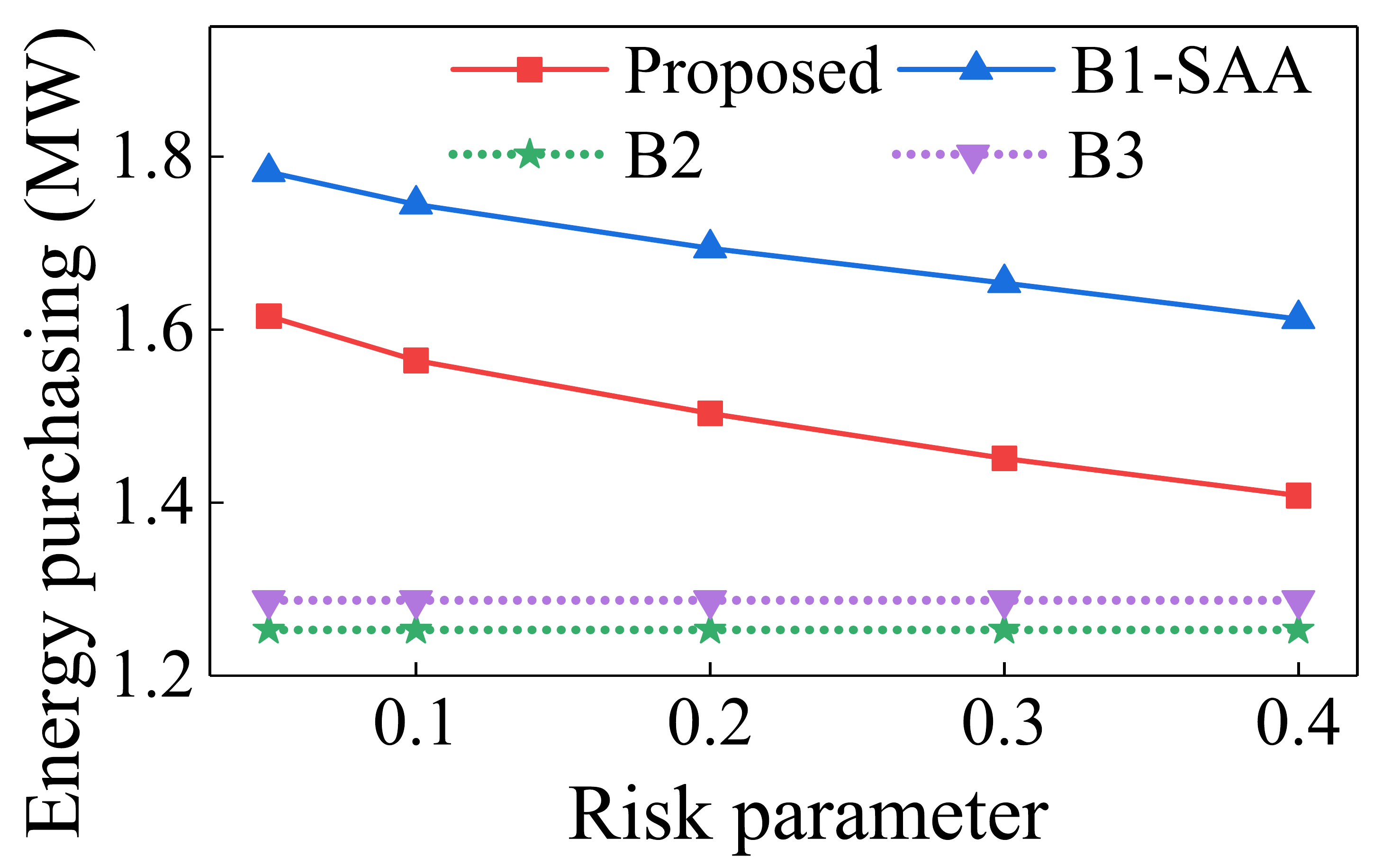}}
	\subfigure[]{\includegraphics[width=0.49\columnwidth]{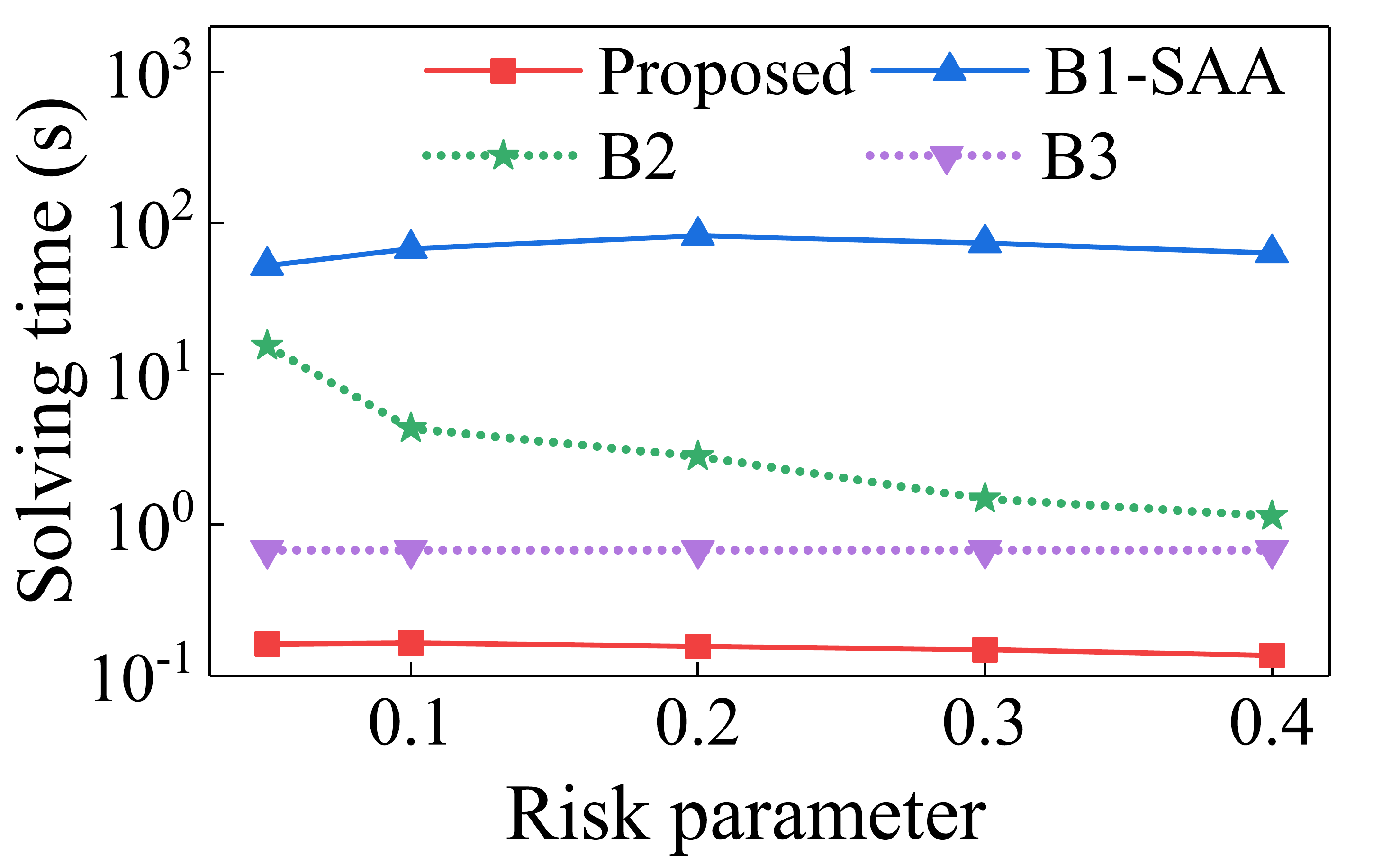}}
	\vspace{-4mm}
 	\caption{Results of (a) average utilization rates of DG, (b) maximum violation probabilities, (c) energy purchasing, and (d) solving times in the case based on the IEEE 123-bus system. Here the uncertainties $\bm \omega$ follow the Weibull distribution.}
	\label{fig_results_123Bus}
	\vspace{-4mm}
\end{figure}

The results of different models on the IEEE 123-bus system are illustrated in Fig. \ref{fig_results_123Bus}. Similar to the results on the 33-bus system, the SOCP relaxation \textbf{B2} can not always guarantee the feasibility of solutions due to the existence of reverse power flows. The risk-neutral model \textbf{B3} also fails to satisfy the JCC because it directly ignores the impacts of uncertainties. For the linearized DistFlow \textbf{B1-SAA}, even though the SAA method is introduced to handle the JCC, its energy purchasing amount is still much higher than that of the proposed model. This result indicates that the linearized DistFlow introduces significant conservativeness and may harm energy efficiency. Moreover, since SAA needs to introduce a large number of binary variables, its computational efficiency is much worse than that of the proposed one. For example, at the risk parameter $\epsilon=0.2$, the solving time of \textbf{B1-SAA} reaches 82.63s, while it is only 0.15s in the proposed model. These results further confirm the great performance of the proposed model.  

\subsection{Sensitivity analysis} \label{sec_sensitive}
In this section, we investigate how the neuron numbers of MLPs affect the performance of the proposed model. The simulations are based on the IEEE 33-bus system.
The hidden layer numbers of MLPs are fixed at three, and the used samples of $\bm \omega$ are the same as those in \textbf{Case 1}.

\subsubsection{Neuron number of quantile-MLP}
\begin{figure}
		\vspace{-4mm}
	\subfigbottomskip=-4pt
	\subfigcapskip=-4pt
	\centering
	\subfigure[]{\includegraphics[width=0.49\columnwidth]{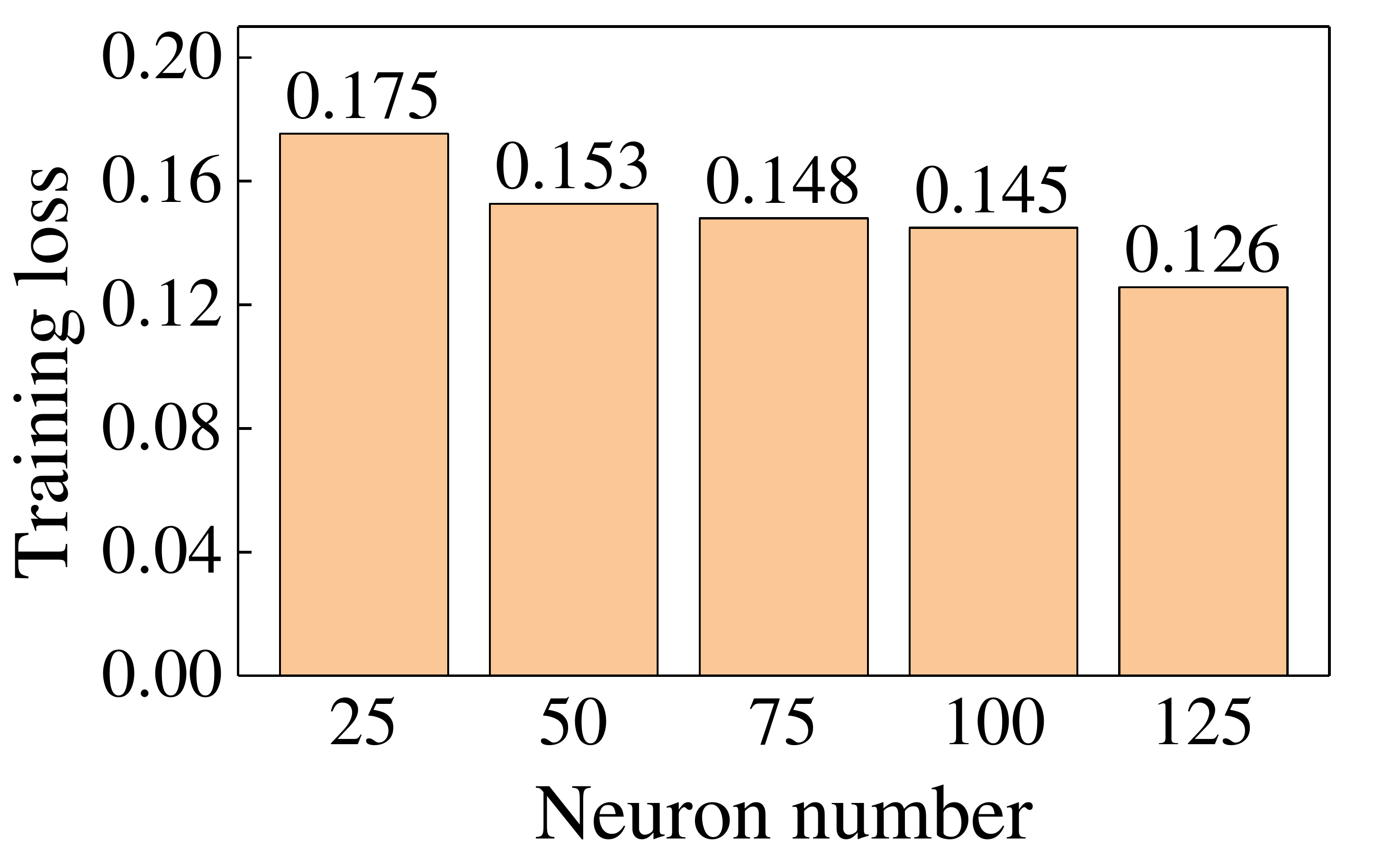}}
	\subfigure[]{\includegraphics[width=0.49\columnwidth]{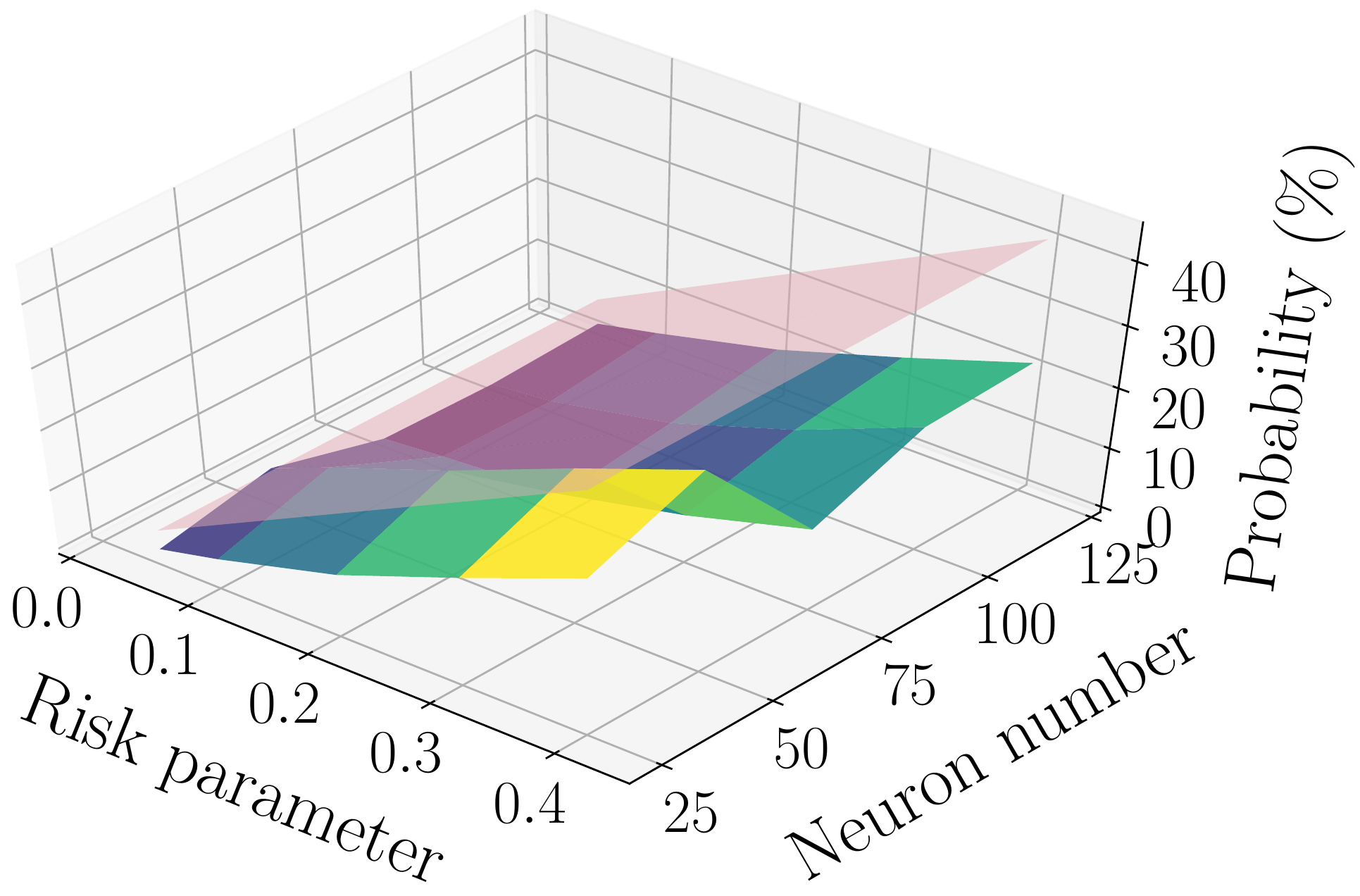}}
	\subfigure[]{\includegraphics[width=0.49\columnwidth]{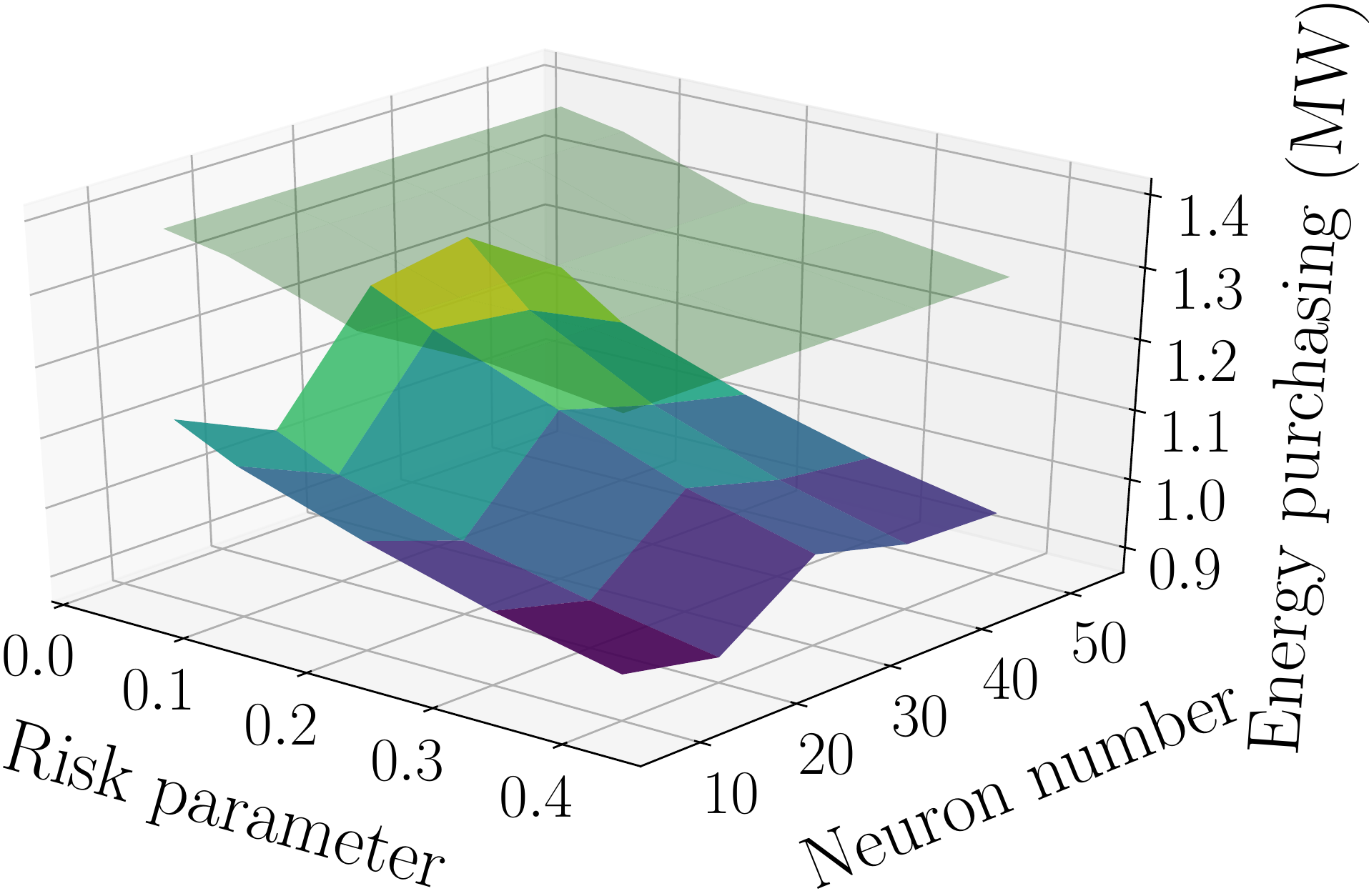}}
	\subfigure[]{\includegraphics[width=0.49\columnwidth]{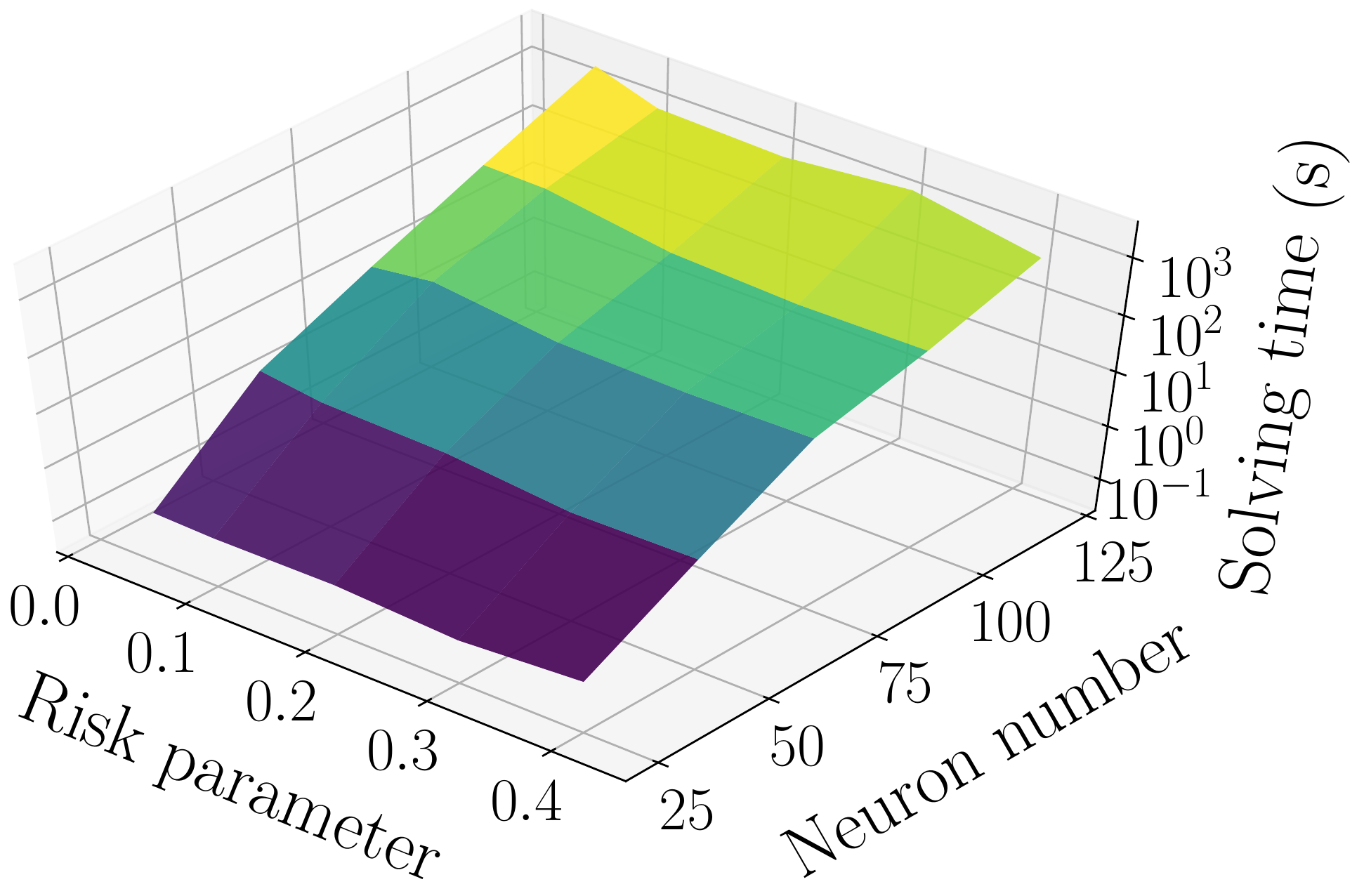}}
	\vspace{-4mm}
 	\caption{Results of (a) loss function of the quantile-MLP, i.e., Eq. (\ref{eqn_loss_quantile}), (b) maximum violation probability, (c) energy purchasing, and (d) solving times under different neuron numbers. In (b), the red surface represents the maximum allowable violation probability, i.e, the given risk parameter. In (c), the blue surface on the top refers to the energy purchasing of \textbf{B1}.}
	\label{fig_results_neuronNum}
	\vspace{-4mm}
\end{figure}

The results of the proposed model with different neuron numbers in the quantile-MLP are illustrated in Fig. \ref{fig_results_neuronNum}, where  ``neuron number" refers to the neuron number in each hidden layer. The structure of the loss-MLP is fixed as (10, 10, 10). With the growth of the neuron number, the approximation ability of the quantile-MLP becomes stronger. Thus, the prediction loss decreases, as shown in Fig. \ref{fig_results_neuronNum}(a). 
Since we use an inner approximation (\ref{eqn_quantile_constraint3}) to replace the original joint chance constraint (\ref{eqn_JCC}) in the calibration step, the maximum violation probabilities are always lower than the risk parameter, i.e., the red surface in Fig. \ref{fig_results_neuronNum}(b). With the growth of the neuron number, the prediction error of the quantile-MLP can be either negative or positive. As a result, both the maximum violation probability and energy purchasing are not monotonous with respect to the neuron number. Nevertheless, the energy purchasing of the proposed model is always lower than that of \textbf{B1}, i.e., the green surface in Fig. \ref{fig_results_neuronNum}(c). 
The solving time grows rapidly with the increase of the neuron number, as illustrated in Fig. \ref{fig_results_neuronNum}(d). According to (\ref{eqn_reformulation}), the integer variable number introduced by reformulating the quantile-MLP is equal to the neuron number. Therefore, a larger neuron number leads to a higher computational burden. Nevertheless, with a small neuron number, the proposed model can already achieve desirable optimality and feasibility simultaneously in a short time, e.g., the solving time is around 0.3s when the neuron number is set as 25.

\subsubsection{Neuron number of loss-MLP}
\begin{figure}
	\subfigbottomskip=-4pt
	\subfigcapskip=-4pt
	\centering
	\subfigure[]{\includegraphics[width=0.49\columnwidth]{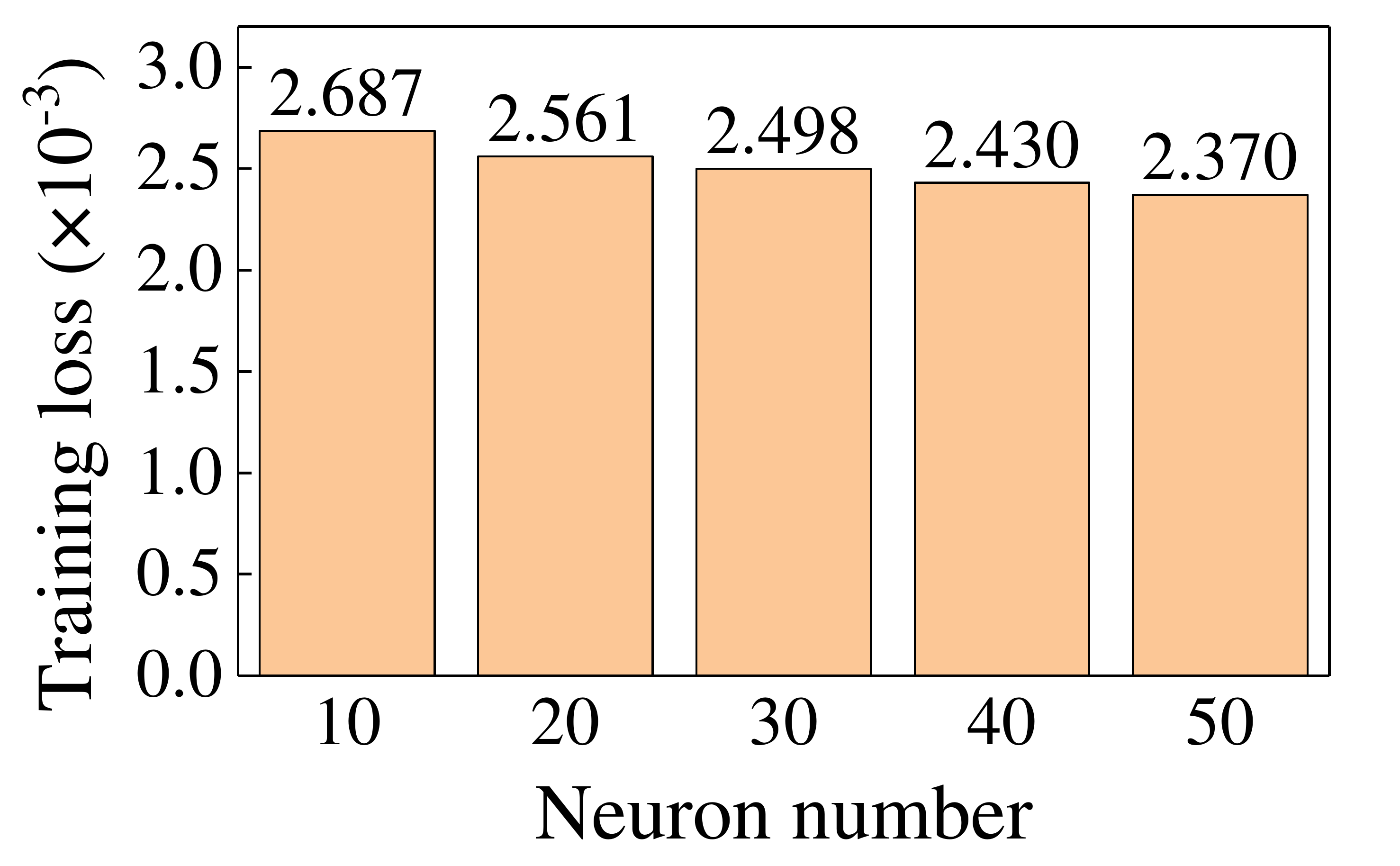}}
	\subfigure[]{\includegraphics[width=0.49\columnwidth]{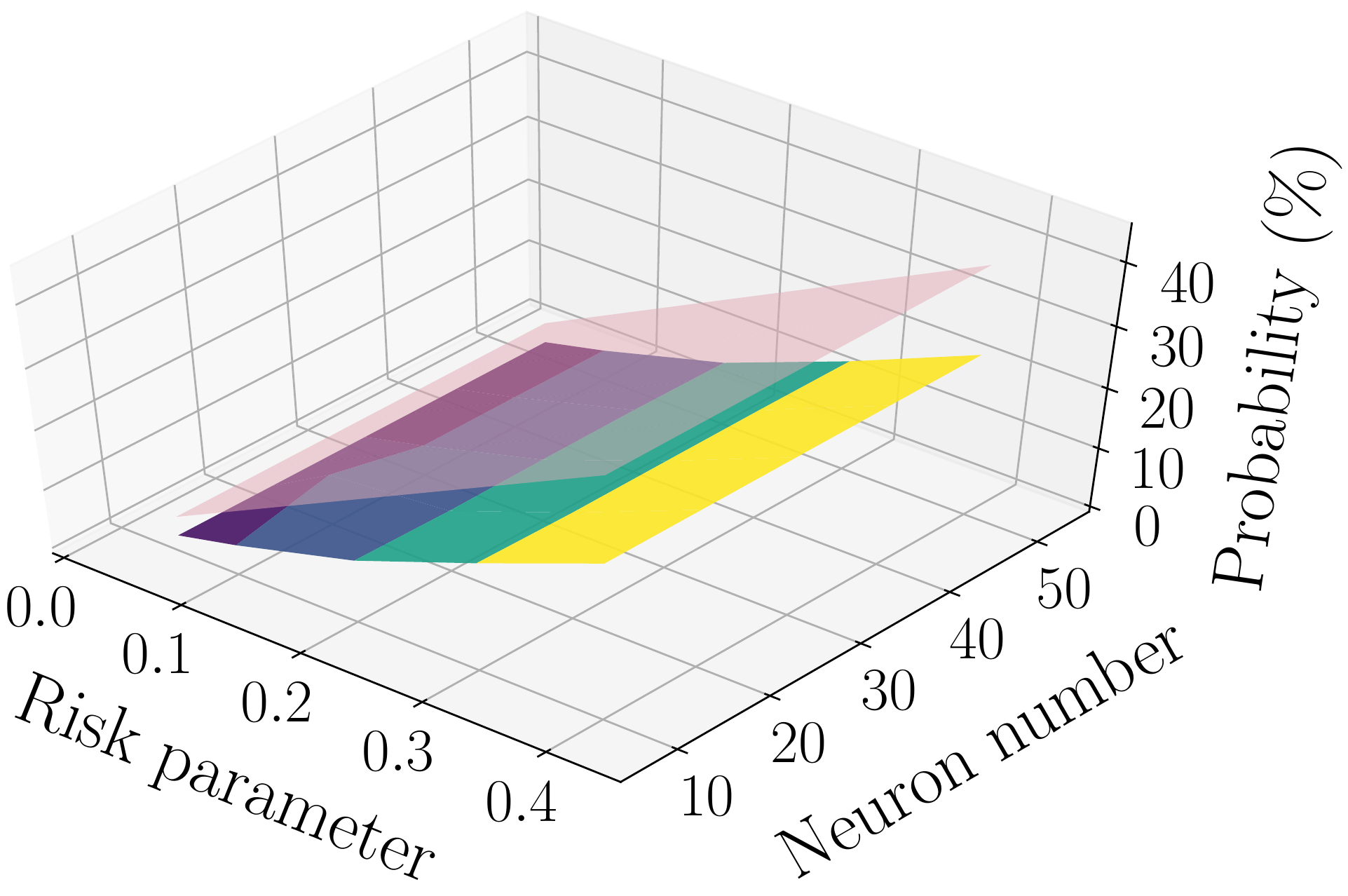}}
	\subfigure[]{\includegraphics[width=0.49\columnwidth]{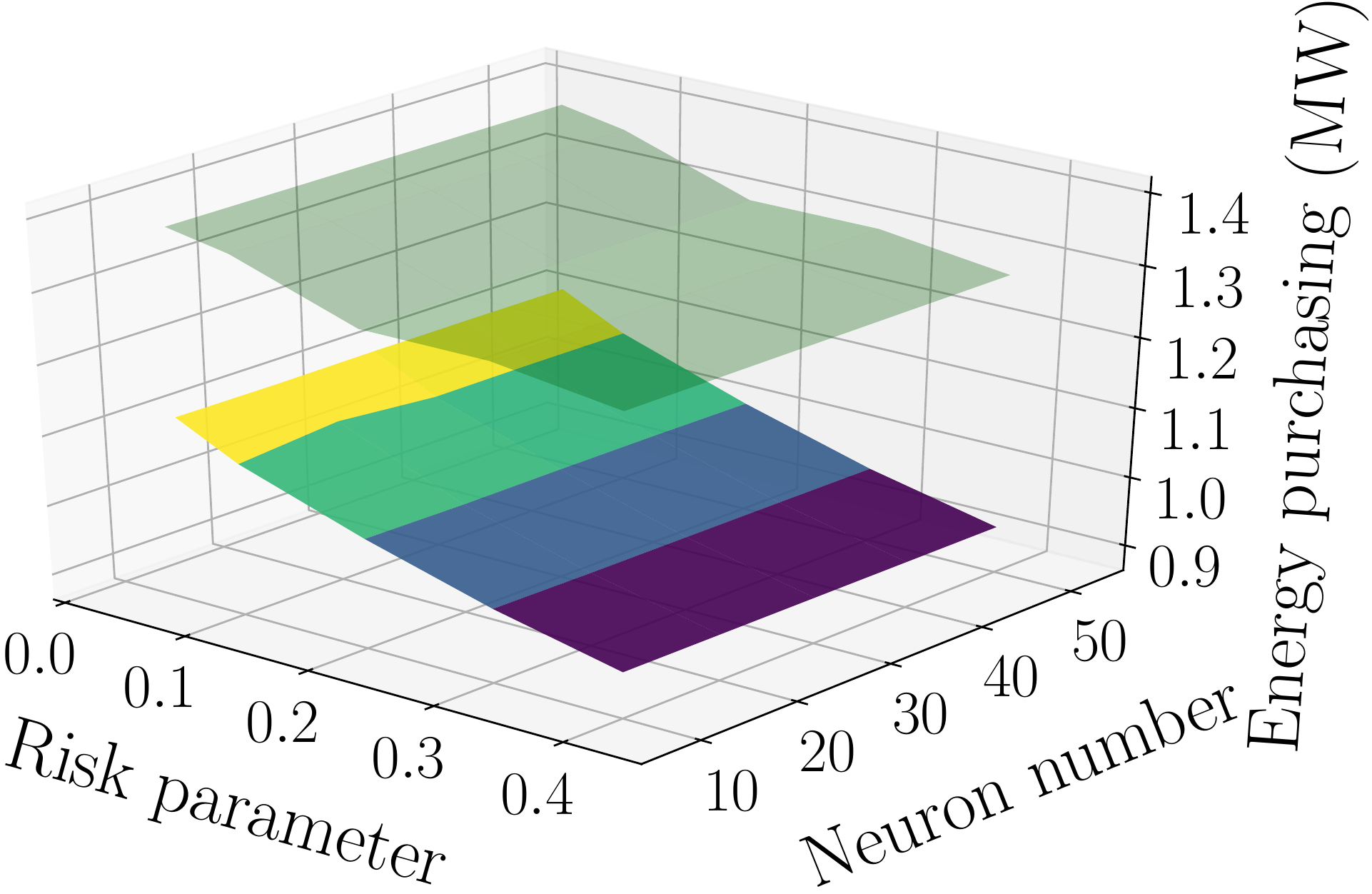}}
	\subfigure[]{\includegraphics[width=0.49\columnwidth]{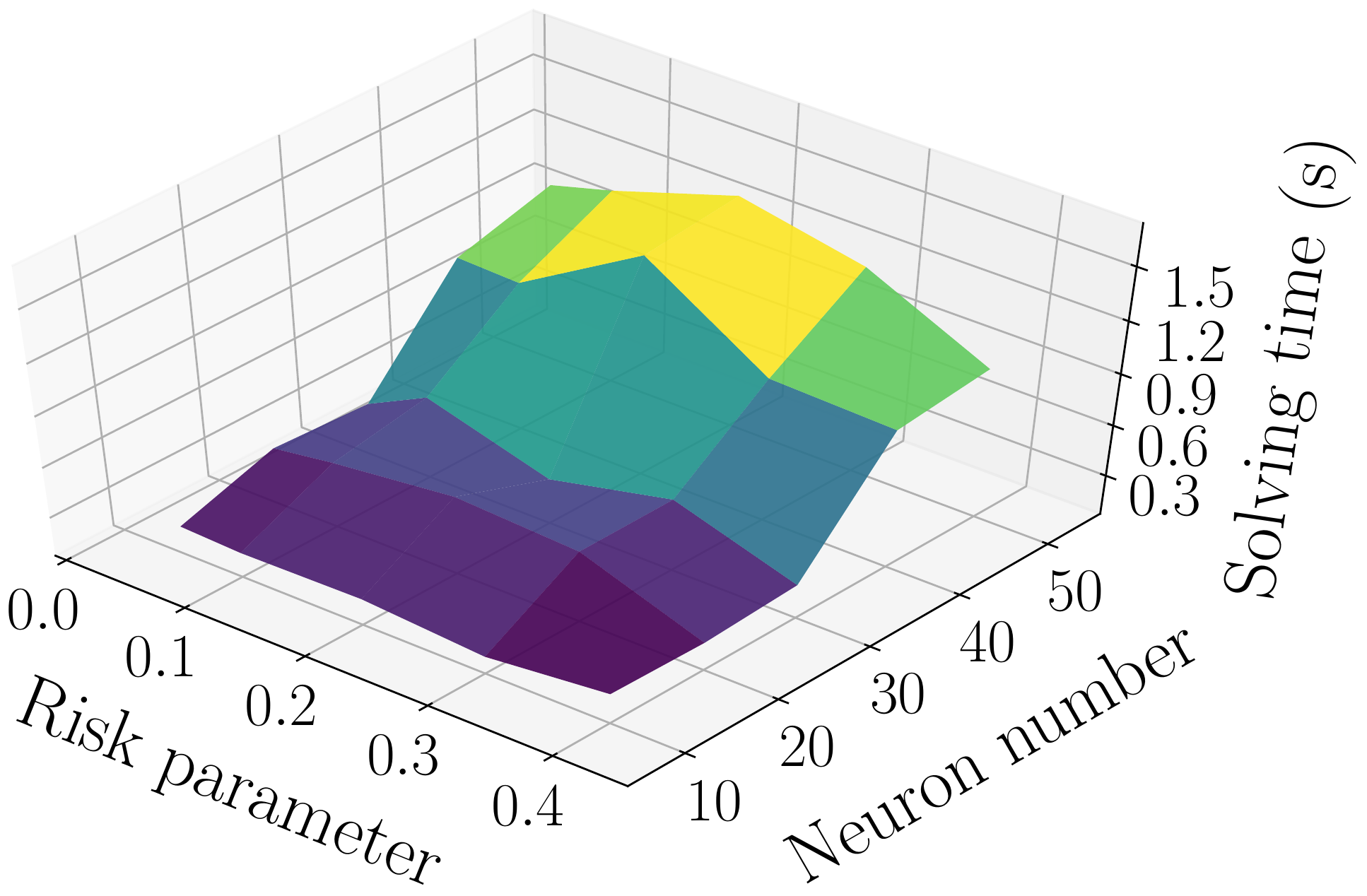}}
	\vspace{-4mm}
 	\caption{Results of (a) loss function of the loss-MLP, (b) maximum violation probability, (c) energy purchasing, and (d) solving times with different neuron numbers in the loss-MLP.}
	\label{fig_results_neuronNum_lossMLP}
	\vspace{-4mm}
\end{figure}

We further investigate the effects of the loss-MLP's neuron number on the proposed model's performance, and the results are summarized in Fig. \ref{fig_results_neuronNum_lossMLP}. Here the neuron number of the quantile-MLP is fixed as (25, 25, 25). Similarly, increasing the neuron number can reduce the training loss of the loss-MLP because this can enhance the prediction accuracy of the loss-MLP, as shown in Fig. \ref{fig_results_neuronNum_lossMLP}(a). However, the power loss is usually much smaller than the summation of power demands. Therefore, even if we change the neuron number, the optimality and feasibility of the proposed model's solutions are nearly constant. Nevertheless, the maximum violation probability is always lower than the required values, i.e., the red surface in Fig. \ref{fig_results_neuronNum_lossMLP}(b), and the energy-efficiency is always better than that of \textbf{B1}, i.e., the green surface in Fig. \ref{fig_results_neuronNum_lossMLP}(c). According to (\ref{eqn_reformulation}), the number of the auxiliary binary variables introduced by reformulating the loss-MLP is equal to its neuron number. Thus, a larger neuron number results in a higher computational burden, and further leads to a longer solving time, as shown in Fig. \ref{fig_results_neuronNum_lossMLP}(d). Nevertheless, a small number of neurons is enough for the proposed model because excellent optimality and feasibility can be already accomplished.

\section{Conclusions} \label{sec_conclusion}
In this paper, we propose a deep-quantile-regression-based surrogate model for the JCC-OPF problem. In the proposed model, two MLPs are trained to predict the $1-\epsilon$ quantile of the maximum constraint violation and expected power loss, respectively. By reformulating the forward propagation of the two MLPs into mixed-integer linear constraints, the JCC-OPF problem can be replicated by the proposed learning-based surrogate model in a mixed-integer form with no need for power network parameters. Two pre-processing steps, i.e., data augmentation and calibration, are further designed to enhance the performance of the proposed model. The data augmentation step trains an XGBoost-based regressor to generator more training samples so that the accuracy of the quantile regression can be improved. The calibration step designs a positive parameter to calibrate the deep quantile regression to improve the feasibility of solutions. Simulation results based on the IEEE 33- and 123-bus distribution systems confirm that the proposed model can successfully replicate the JCC-OPF problem without the network parameters. Moreover, its optimality is better than the widely used linearized DistFlow model under arbitrary uncertainties, while its feasibility performance is much greater than the SOCP relaxation of AC OPF. Numerical experiments also demonstrate that a small number of neurons is already enough for the proposed model to achieve optimality and feasibility, so computational efficiency can be also guaranteed.

\appendices
\setcounter{table}{0}   
\renewcommand\thetable{\Alph{section}\arabic{table}}  

\section{} \label{app_1}
\emph{Proof of \textbf{Proposition} \ref{proposition_1}}: 
The term on the right-hand side of (\ref{eqn_proposition_1}) is equal to
\begin{align}
\mathbb E (\text{Loss}^\text{QR}) =  -\epsilon\int_{-\infty}^{\hat{\mathcal{Q}}^{1-\epsilon}}(h - \hat{\mathcal{Q}}^{1-\epsilon})dF_{H}(h)  \notag \\
+ (1-\epsilon)\int_{\hat{\mathcal{Q}}^{1-\epsilon}}^{\infty}(h - \hat{\mathcal{Q}}^{1-\epsilon})dF_{H}(h), \label{eqn_proposition_2}
\end{align}
where $F_{H}(\cdot)$ denotes the cumulative distribution function of $h(\bm x, \bm \omega)$ at $\bm x$ under the uncertainty $\bm \omega$. At the optimal solution that minimizing the expectation (\ref{eqn_proposition_2}), the derivative of the expectation loss should be zero:
\begin{align}
\left. \frac{\partial \mathbb E (\text{Loss}^\text{QR})}{\partial \hat{\mathcal{Q}}^{1-\epsilon}}\right|_{y}=0, \label{eqn_proposition_3}
\end{align}
where $y$ is the optimal solution of $\hat{\mathcal{Q}}^{1-\epsilon}$. Then, 
by substituting (\ref{eqn_proposition_2}), Eq. (\ref{eqn_proposition_3}) can be converted into the following form based on the Leibniz integral rule:
\begin{align}
&\epsilon \int_{-\infty}^{y} dF_{H}(h) - (1-\epsilon)\int_{y}^{\infty} dF_{H}(h)=0. \label{eqn_proposition_4}
\end{align}
By substituting $F_{H}(-\infty)=0$ and $F_{H}(\infty)=1$, Eq. (\ref{eqn_proposition_4}) can be further reformulated as:
\begin{align}
F_{H}(y) = 1 - \epsilon \Leftrightarrow y = \mathcal{Q}_{\bm \omega}^{1-\epsilon}(h(\bm x, \bm \omega)).
\end{align}
This completes the proof.

\section{} \label{app_2}
\emph{Proof of \textbf{Proposition} \ref{proposition_2}}: 
Based on (\ref{eqn_loss_powerLoss}), the expectation of $\text{Loss}^\text{pl}$ can be expressed as:
\begin{align}
\mathbb E (\text{Loss}^{\text{pl}}) &= \mathbb E \left((p^\text{loss} - \hat{p}^\text{loss}(\bm x))^2\right)\notag\\
&= \mathbb E \left((p^\text{loss}-\mathbb{E}_{\bm \omega}(p^\text{loss}))^2\right) + \left(\mathbb{E}_{\bm \omega}(p^\text{loss}) - \hat{p}^\text{loss}(\bm x)\right)^2 \notag\\
&=\text{Var}(p^\text{loss}) + \left(\mathbb{E}_{\bm \omega}(p^\text{loss}) - \hat{p}^\text{loss}(\bm x)\right)^2, \label{eqn_A2_1}
\end{align} 
where $\text{Var}(p^\text{loss})$ is the variance of $p^\text{loss}$. By regarding $\mathbb E (\text{Loss}^{\text{pl}})$ as a function of $\hat{p}^\text{loss}(\bm x)$, the minimum value of $\mathbb E (\text{Loss}^{\text{pl}})$ occurs at $\hat{p}^\text{loss}(\bm x)=\mathbb{E}_{\bm \omega}(p^\text{loss})$ according to (\ref{eqn_A2_1}). This completes the proof.

\footnotesize
\bibliographystyle{ieeetr}
\bibliography{ref}
\end{document}